\documentclass[a4paper,10pt]{article}
\usepackage[dvipsnames]{xcolor}

\usepackage{amssymb}
\usepackage{graphicx}
\usepackage{bbm}
\usepackage{mathrsfs}
\usepackage{amsmath}
\usepackage{amsthm}
\usepackage[center]{caption}
\usepackage{enumerate}
\usepackage{verbatim}
\usepackage{authblk}
\usepackage[hmargin=2.5cm, vmargin=3.5cm]{geometry}

\usepackage{bm}

\graphicspath{ {C:/Users/sgj22/Pictures/latexgraphs/FDD/} }
\usepackage{color}
\usepackage{tikz}
\usepackage[normalem]{ulem}
\usepackage{hyperref}
\hypersetup{
    colorlinks,
    linkcolor={red!59!black},
    citecolor={blue!59!black},
    urlcolor={blue!89!black}
}

\numberwithin{equation}{section}

\newtheorem{thm}{Theorem}[section]
\newtheorem{cor}[thm]{Corollary}
\newtheorem{lem}[thm]{Lemma}
\newtheorem{prop}[thm]{Proposition}
\newtheorem{conj}[thm]{Conjecture}

\newtheorem{defn}[thm]{Definition}
\newtheorem{rmk}[thm]{Remark}

\newcommand{\BeG}{\begin{align*}}
\newcommand{\EnD}{\end{align*}}
\newcommand{\bae}{\begin{align}}
\newcommand{\eae}{\end{align}}

\renewcommand{\P}{\mathbb{P}}
\newcommand{\E}{\mathbb{E}}

\newcommand{\ind}{\mathbbm{1}}

\newcommand{\bp}{\begin{proof}}
\newcommand{\ep}{\end{proof}}
\def\bal#1\eal{\begin{align*}#1\end{align*}}

\newcommand{\Z}{\mathbb{Z}}

\renewcommand{\vec}[1]{\mathbf{#1}}

\newcommand{\vx}{\vec{x}}
\newcommand{\vy}{\vec{y}}
\newcommand{\vz}{\vec{z}}
\newcommand{\uk}{\uparrow k}
\newcommand{\Uk}{\Uparrow k}

\newcommand{\lf}{\lfloor}
\newcommand{\rf}{\rfloor}

\DeclareMathOperator{\Tr}{Tr}
\DeclareMathOperator{\Leb}{Leb}
\DeclareMathOperator{\Vol}{Vol}
\DeclareMathOperator{\wind}{wind}
\DeclareMathOperator{\Realpart}{Re}

\author{\textsc{Samuel G.G.~Johnston and Neil O'Connell}}

\title{Scaling limits for non-intersecting polymers and Whittaker measures}

\begin{document}
\maketitle

\begin{abstract}
We study the partition functions associated with non-intersecting polymers in a random environment. 
By considering paths in series and in parallel, the partition functions carry natural notions of subadditivity, allowing the effective study of their asymptotics.
For a certain choice of random environment, the geometric RSK correspondence provides an explicit representation of the partition functions in terms of a stochastic interface.  Formally this leads to a variational description of the macroscopic behaviour of the interface and hence the free energy of the associated non-intersecting polymer model.  
At zero temperature we relate this variational description to the Mar\v{c}enko-Pastur distribution, and give a new derivation of the surface
tension of the bead model.
\end{abstract}


\section{Introduction and summary} \label{sec int}

We study the partition functions associated with a natural model for non-intersecting polymers in a random environment. 
Apart from being an interesting physical model in its own right, this model is motivated by recent developments on
connections between random polymers and Whittaker functions, obtained via the geometric RSK correspondence~\cite{COSZ,OSZ}.
By considering paths in series and in parallel, the partition functions carry natural notions of subadditivity allowing the 
effective study of their asymptotics.  We use this subadditivity to show that for a small number of paths, the free energy has a linear
dependence on the number of paths (Theorem \ref{finite k thm 2}).  Interestingly, it seems quite difficult to prove this using 
the integrable structure of the log-gamma polymer model; instead we give a general argument, valid for general weights.
We also determine the correct order of scaling for the partition function associated with a large number of paths (see \eqref{LL}, \eqref{UL} and Corollary \ref{correct scaling}). 

We then turn our discussion to the random polymer model with log-gamma weights, where the geometric RSK correspondence provides an explicit representation 
of the partition functions in terms of a stochastic interface \cite{COSZ,fun}.  Formally this leads to a variational description of the macroscopic behaviour of the interface and hence the free energy of the associated non-intersecting polymer model.  
At zero temperature we relate this variational description to the Mar\v{c}enko-Pastur distribution, and give a new derivation of the surface
tension of the bead model, which was recently computed by Sun \cite{sun} using very different methods.
In the remainder of Section 1, we summarise the main results of the paper, beginning with some preliminary definitions.

\subsection{Preliminary definitions}
Let $\mathbf{e}_1 = (1,0)$ and $\mathbf{e}_2 = (0,1)$ be the standard basis for $\mathbb{Z}^2$. For points $x$ and $y$ in $\mathbb{Z}^2$, a path $\pi$ from $x$ to $y$ is a set of points $\pi = \{ a_0,\ldots,a_p \}$ in $\Z^2$ such that $a_0 = x, a_p = y$, and $a_{j+1} - a_j$ is equal to either $\mathbf{e}_1$ or $\mathbf{e}_2$ for every $1 \leq j < p$.  

We say a vector $\vec{x} = (x_1,\ldots,x_k) \in \mathbb{Z}^{2 \times k}$ is a $k$-point if $x_1,\ldots,x_k$ are distinct points in $\mathbb{Z}^2$. Suppose $\vx = (x_1,\ldots,x_k)$ and $\vy = (y_1,\ldots,y_k)$ are $k$-points and for each $i$, $\pi_i$ is a path from $x_i$ to $y_i$. If the sets $\pi_1,\ldots,\pi_k$ are disjoint, we say the $k$-tuple $\pi = (\pi_1,\ldots,\pi_k)$ is non-intersecting and refer to $\pi$ as a $k$-path. We write $\Gamma_{\vx \to \vy}$ for the set of $k$-paths from $\vx $ to $\vy$, and write $\vx \leq \vy$ whenever $\Gamma_{\vx \to \vy}$ is non-empty. (We emphasise that a $k$-path will always refer to a non-intersecting $k$-tuple of paths travelling north-east in $\mathbb{Z}^2$.) 

Suppose further we have a random environment $\{ \omega(z) : z \in \Z^2 \}$---a collection of independent and identically distributed finite-expectation random variables under a probability measure $\P$. For a $k$-path $\pi = (\pi_1,\ldots,\pi_k)$ from $\vx$ to $\vy$, we define the energy $F(\pi)$ of $\pi$ to be the sum of the weights included in the $k$-path, not including the starting points. That is, $F(\pi)$ is the random variable
\begin{align} \label{FE}
F(\pi) := \sum_{i=1}^k \sum_{ z \in \pi_i - \{ x_i \} } \omega(z).
\end{align}
Suppose $\beta$ is a real number such that $G(\beta) := \E [ e^{\beta \omega(0,0)} ] < \infty$.
The goal of the present paper is to study the partition function
\begin{align} \label{PF}
Z_{\vec{x} \to \vec{y}}(\beta) := \sum_{ \pi \in \Gamma_{\vx \to \vy}  } \exp \left( \beta F(\pi) \right)
\end{align}
associated with the non-intersecting directed random polymer running from $\vx$ to $\vy$, with a particular emphasis on the asymptotics of $Z_{\vx \to \vy}(\beta)$ as the lengths of the paths and/or the number of paths grows to infinity. We emphasise that the assumptions $\nu := \E[\omega(0,0)] < \infty$ and $G(\beta) = \E [ e^{\beta \omega(0,0)} ]  < \infty$ are in force throughout the paper.

This article is divided into three interconnected parts. The first part is dedicated to proving asymptotic results for the partition functions in a general environment. In the second part we introduce a collection of random functions on the $N \times N$ square known as stochastic interfaces, and express the non-intersecting partition functions associated with a particular random environment in terms of a stochastic interface, developing connections with random matrix theory in the process. In the final part, we study the scaling limits of these interfaces as $N \to \infty$, relating our model to the Mar\v{c}enko-Pastur law and Gaussian fields, and developing the machinery along the way to obtain new derivation of the surface tension of the bead model. In the remainder of Section \ref{sec int} we overview each of these three respective parts in more detail.

Before proceeding, we mention a couple of conventions we use to lighten notation. For real numbers $u$, let $\lfloor u \rfloor$ denote the largest integer less than or equal to $u$, and let $\lceil u \rceil $ denote the smallest integer greater than or equal to $u$. Whenever $(a_n)_{n \in \mathbb{Z}}$ are variables indexed by the integers, for real numbers $r$ we set $a_r := a_{\lfloor r \rfloor}$, where $\lfloor r \rfloor$ is the largest integer not greater than $r$. We will use similar conventions throughout, in particular
\begin{align*}
\text{whenever $c$ is positive and $N$ is an integer, $(N,cN)$ always refers to the point $(N, \lfloor cN \rfloor)$ of $\mathbb{Z}^2$. }
\end{align*}
If $f$ is a function defined on a subset of $\mathbb{Z}^2$ or $\mathbb{R}^2$, we say $f$ is symmetric if $f(r,s) = f(s,r)$ for every $s$ and $r$. Finally, if $\vx$ is a $k$-point and $a$ is an element of $\mathbb{Z}^2$, define the translation $\vx + a$ to be the $k$-point whose $i^{\text{th}}$ component is given by $\vx_i + a$.

\subsection{Non-intersecting polymers} \label{polymer o}
The first part of this article is concerned with studying key properties of non-intersecting polymers in a general environment, with a particular focus on developing tools to study their asymptotics. These results are stated in full in Section \ref{polymer}, though we provide a brief outline here.

We say a $k$-point $\vx = (x_1,\ldots,x_k)$ is nice if either $x_{i+1}$ lies strictly north and strictly west of $x_i$, or $x_{i+1} = x_i + \mathbf{e}_2$. This is a technical condition which roughly speaking ensures that $k$-paths going to and from a $k$-point $\vx$ do not intersect. (See the beginning of Section \ref{polymer} for a precise statement and proof.)

Our main tools for tackling the asymptotics of non-intersecting polymer partition functions $Z_{\vx \to \vy}(\beta)$ are the dual notions of series and parallel concatenation, which we now outline. Given a $k$-path $\pi$ from $\vx$ to $\vy$ and $k$-path $\gamma$ from $\vy$ to $\vz$, we may concatenate the paths $\pi$ and $\gamma$ in series to form a new $k$-tuple $\pi \oplus \gamma$ of paths from from $\vx$ to $\vz$. Provided the intermediate $k$-point $\vy$ is nice, the concatenated path $\pi \oplus \gamma$ turns out to also be non-intersecting, and in this case the operation $\oplus$ is additive in the sense that  
\begin{align*}
F( \pi \oplus \gamma) = F( \pi) + F( \gamma),
\end{align*}
where $F$ is given as in \eqref{FE}.
In particular, we have an $F$-preserving injection taking every pair of paths in $\Gamma_{\vx \to \vy } \times \Gamma_{\vy \to \vz}$ to a new path $\pi \oplus \gamma$ in $ \Gamma_{\vx \to \vz}$. In Section \ref{polymer}, we use this injection to prove the series bound, which states that for nice $k$-points $\vx \leq \vy \leq \vec{z}$, the partition functions satsify
\begin{align} \label{series bound}
Z_{\vx \to \vz}(\beta) \geq Z_{\vx \to \vy}(\beta) Z_{\vy \to \vz}(\beta).
\end{align}
The series bound gives rise to subadditivity in the logarithmic partition functions $\log Z_{\vx \to \vy}(\beta)$, which we use in conjunction with Kingman's subadditive ergodic theorem \cite{kin} to obtain the logarithmic asymptotics of $Z_{\vx \to \vy}(\beta)$ as $\vx$ and $\vy$ grow far apart in the $(1,c)$ direction. Theorem \ref{finite k thm} states that for any pair of nice $k$-points $\vx$ and $\vy$ 
\begin{align} \label{finite k ex}
\frac{1}{N} \log Z_{\vec{x} \to \vec{y} + (N , c N ) }(\beta) \text{ converges to a deterministic limit $f_c(k,\beta)$},
\end{align}
where the constant $f_c(k,\beta)$ is independent of $\vx$ and $\vy$. As we will see later, it is a natural consequence of the parallel bound that the limit $f_c(k,\beta)$ is concave in the $k$-variable, in that for any $j,k$ we have
\begin{align*}
f_c(k+j,\beta) \leq f_c(j,\beta) + f_c(k,\beta).
\end{align*}
However, it turns out that $f_c(k,\beta)$ has a linear dependence on $k$. Indeed, Theorem \ref{finite k thm 2} states that 
\begin{align} \label{finite k ex}
f_c(k,\beta) = k  f_c(\beta),
\end{align}
where we are writing $f_c(\beta)$ for $f_c(1,\beta)$. 

Where series concatenation involves combining paths by glueing their ends together, parallel concatenation combines paths by letting them run side-by-side. Namely, any $(k+j)$-path $\pi$ has a unique decomposition $\pi = \pi' \boxplus \pi''$ where $\pi'$ is a $k$-path and $\pi''$ is a $j$-path. Again this procedure is additive in that $F(\pi' \boxplus \pi'') = F(\pi') + F(\pi'')$, and in Section \ref{polymer 1} we use it prove the parallel bound 
\begin{align} \label{parallel bound}
Z_{\vx \to \vy}(\beta) \geq Z_{\vx' \to \vy'}(\beta) Z_{\vx'' \to \vy''}(\beta),
\end{align}
where $\vx := (x_1,\ldots,x_{k+j})$, $\vx' = (x_1,\ldots,x_k)$ and $\vx'' = (x_{k+1},\ldots,x_{k+j})$, and $\vy, \vy'$ and $\vy''$ are defined similarly. 

The parallel bound \eqref{parallel bound} gives a subadditivity in the $k$- variable which allows us to prove a result analogous to \eqref{finite k ex} concerning infinitely many paths of finite length running side-by-side. Thereafter, we study the infinite temperature $(\beta = 0)$ case, where the Lindstr{\"o}m-Gessel-Viennot formula may be used in conjunction with Szeg{\"o}'s powerful limit theorem to obtain explicit expressions for the large-$k$ logarithmic asymptotics of the partition funcitions in terms of a Laurent series with binomial coefficients.

Finally we look at the case of many long paths. By studying the large-$N$ asymptotics of a variant of Macmahon's formula \cite{sta}, in Lemma \ref{infinite growth} we obtain explicit logarithmic asymptotics for the partition functions associated with the infinite temperature case. We then use Jensen's inequality to relate the positive temperature partition functions with their infinite temperature counterparts (Lemma \ref{infinite bounds}), leading us to Corollary \ref{correct scaling}, which establishes $N^2$ behaviour for the logarithmic partition function when both the number of paths aswell as their lengths are of order $N$. 

\subsection{Stochastic interfaces and Whittaker measures} \label{SIM o}
In the second part of this article---the results of which we state in full in Section \ref{SIM}---we introduce stochastic interfaces, random functions defined on the square $S_N := \{1,\ldots,N\}^2$ of $\mathbb{Z}^2$. With the exception of Theorem \ref{large mu} which is new, our work in this section is dedicated to reformulating results from Baryshnikov \cite{bar} and Corwin, O'Connell, Sepp{\"a}l{\"a}inen and Zygouras \cite{COSZ} in the language of stochastic interface models.

For a convex interaction potential $V$, the stochastic interfaces we consider are essentially random functions $\phi^N:S_N \to \mathbb{R}$ with distributions proportional to
\begin{align} \label{D law}
\exp \left( -  \sum_{ \langle x,y \rangle }V( \phi(y) - \phi(x) ) - \mu \sum_{ x \in D_N } \phi(x)  \right),
\end{align}
where $D_N := \{ (i,i) : i = 1, \ldots, N \}$ is the diagonal of $S_N$, and $\langle x,y \rangle$ is the set of edges in $S_N$ directed in the north or east direction. The central idea connecting directed polymers with stochastic interfaces is Proposition \ref{interface 2}, a restatement of a result by Corwin et al \cite{COSZ}, which states that the non-intersecting partition functions associated with a certain random environment may be jointly expressed in terms of a stochastic interface. 

More specifically, let $\tau_\mu^N(m,k)$ be the partition function associated with a $k$-paths in an $m \times N$ rectangle, where the variables in the random environment are log-gamma distributed in that each $e^{\beta \omega(z)}$ has the inverse gamma law
\begin{align} \label{lg1}
I_\mu(ds) = \frac{ s^{ - \mu - 1} e^{ - 1/s} }{ \Gamma(\mu) } \ind_{s > 0} ds
\end{align}
of parameter $\mu > 0$. We now define a random function $\varphi_\mu^N:S_N \to \mathbb{R}$ in terms of the joint partition functions $\tau_\mu^N(m,k)$ by setting 
\begin{align*}
\varphi_\mu^N(i,j) := \log \left( \frac{ \tau^N_\mu({N-j+i,i})}{ \tau^N_\mu({N-j+i,i-1}) }   \right)
\end{align*}
for $i \leq j$ (with a similar definition holding for $i \geq j$). Based on the main result in \cite{COSZ}, we obtain the statement of Proposition \ref{interface 2}: 
\begin{align*}
\varphi^N_\mu \text{ is a stochastic interface with interaction potential $V(u) = \exp(u)$}.
\end{align*}
The small-$\mu$ and large-$\mu$ deformations of the stochastic interface $\varphi^N_\mu$ have respective connections with the eigenvalues ensembles of random matrix theory and the partition functions of deterministic non-intersecting directed polymers. In the former case, by appealing to ideas in Baryshnikov \cite{bar}, we show that an eigenvalue process $\varphi^N_{\mathsf{LUE}}$ associated with the minors of a random matrix from the Laguerre Unitary Ensemble may be expressed in terms of a stochastic interface with the hard-core interaction potential $V(u) := \infty \ind_{ u > 0}$. We use the results of \cite{COSZ} to obtain Proposition \ref{small mu}, which states that
\begin{align} \label{small mu statement}
\mu \varphi^N_\mu \text{ converges in distribution to $\varphi^N_{\mathsf{LUE}}$  as $\mu \to 0$}.
\end{align}
As is mentioned in \cite{COSZ}, by using a tropicalisation of the underlying polymer, one can use the statement \eqref{small mu statement} to recover a result by Johannson \cite{joh} relating non-intersecting last passage percolation with exponential weights to the eigenvalues of a Laguerre matrix.

We have contrasting behaviour as $\mu \to \infty$. Namely, if we define the tilted interface by $\theta^N_\mu(i,j ) := \varphi_\mu^N(i,j)  + (2N+1 - (j+i) ) \log \mu$, then according to Theorem \ref{large mu} we have 
\begin{align} \label{large mu statement}
\theta^N_\mu \text{ converges almost surely to a deterministic function $\theta_{\mathsf{min}}^N$ as $\mu \to \infty$}.
\end{align}
Moreover, the deterministic limit $\theta_{\mathsf{min}}^N$ is the minimiser of a discrete variational problem on the square $S_N$, and has an explicit representation in terms of factorials. 

\subsection{Scaling limits for stochastic interfaces} 
 \label{thermodynamics o}
The final part of this article, which is discussed in full in Section \ref{thermodynamics}, is concerned with using variational heuristics to study the large-$N$ asymptotics of stochastic interfaces. 
In the interest of maintaining a steady flow of ideas, in this part of the article we will be content to provide a plausibility argument based on physical heuristics in place of a rigorous mathematical proof. 

Following the exposition in  Funaki and Spohn \cite{FS}, our first tool for studying the macroscopics of stochastic interfaces is that of surface tension---an asymptotic measure of the energy cost for an interface to lie a certain tilt. With this concept at hand, we conjecture that for certain choices of interaction potential, the macroscopic shapes of certain square interfaces may be given as solutions to variational problems on the unit square $S := [0,1]^2$. In particular we anticipate that if $\bar{\varphi}_\mu^N:S \to \mathbb{R}$ is a rescaling of the interface associated with the partition functions of the non-intersecting log-gamma polymer---with both the square $S_N$ and the height of the interface both rescaled by $\frac{1}{N}$---then  
\begin{align*}
\bar{\varphi}_\mu^N \text{ converges to a deterministic limit $\xi_\mu:S \to \mathbb{R}$ as $N \to \infty$}.
\end{align*}
Moreover, if we let $\sigma^{\mathsf{exp}}$ be the surface tension associated with the exponential interaction potential $V(u) = \mathsf{exp}(u)$, then we expect that the deterministic limit $\xi_\mu$ is the the minimiser of 
\begin{align*}
\text{ $\mathcal{E}_\mu[v] :=  \int_S \sigma^{\mathsf{exp}}( \nabla v(s,t) ) ds dt + \mu \int_0^1 v(s,s) ds,$ }
\end{align*} 
where the minimisation is taken over all functions $v:S \to \mathbb{R}$ satisfying $v(1,1) = 0$. 
Finally, the limit shape is related to the asymptotic partition functions of the non-intersecting directed log-gamma polymer through the following equation. Let $\tau^N_\mu(m,k)$ be the partition function associated with $k$ paths on an $N \times m$ rectangle. Then these arguments suggest that
\begin{align} \label{many path 2}
\lim_{N \to \infty} \frac{1}{N^2} \log \tau^N_\mu( cN, \alpha N) := \int_0^\alpha \xi_\mu(u, 1- c  +u) du. 
\end{align}
Having taken $N \to \infty$, we now consider the behaviour of $\xi_\mu$ as the parameter $\mu$ varies. On the one hand, by considering the relationship the interface $\xi_\mu^N$ has with the eigenvalues of Laguerre matrices, taking a large-$N$ analogue of \eqref{small mu statement} we anticipate that 
\begin{align*}
\text{ $\mu \xi_\mu$  converges to  $\xi_{\mathsf{mp}}$ as $\mu \to 0$},
\end{align*}
where $\xi_{\mathsf{mp}}$ has an explicit expression in terms of the Mar\v{c}enko-Pastur distribution. On the other hand, in the large-$\mu$ case we sketch an argument suggesting that 
\begin{align} \label{psi conv 2}
 \text{$\xi_\mu (s,t)  + (2 - s - t) \log \mu $  converges to  $\xi_{\mathsf{ht}}(s,t)$ as $\mu \to \infty$}, 
\end{align}
where $\xi_{\mathsf{ht}}(s,t)$ has an explicit expression in terms of the function $q(u) = u \log  u$. Given their own direct representations as scaling limits of stochastic interfaces, we also expect both limit functions $\xi_{\mathsf{mp}}$ and $\xi_{\mathsf{ht}}$ to be solutions of explicit variational problems.

Thereafter, we take a more refined argument to study the high-temperature limit --- using the central limit theorem in place of the law of large numbers --- with $\mu$ and $N$ sent to $\infty$ together through the scaling limit $
\mu = \kappa N^2.$
Indeed, Conjecture \ref{conj gauss} states that the random processes $H^N_\kappa := \{ H^N_\kappa(s,t) : 0 \leq s \leq t \leq 1 \}$ given by
\begin{align*} 
H_\kappa^N(s,t) := k(N+m-k) \left( \log \kappa + 2 \log N  \right) + \log \tau_{\kappa N^2}^N(m,k) - \log \# \Gamma^N(m,k)
\end{align*} 
converge in distribution to a centred Gaussian process, and gives a prediction for the covariance functions in terms of the asymptotic densities of a path model. 

Finally, having developed the surface tension framework to study the macroscopics of stochastic interfaces, we use a version of the classical semicircle law of random matrix theory to reverse-engineer a straightforward derivation of a formula for the surface tension associated with the bead model \cite{bou}. To our knowledge, the only other place this formula appears in the literature is in the work of Sun \cite{sun}, who provides a sophisticated derivation involving scaling limits of dimer models \cite{CKP}.

\subsection{Outline of the paper}
The remainder of the paper is structured as follows. In Section \ref{polymer}, we discuss the asymptotics of non-intersecting directed polymers in random environments with a general weight distribution, giving full statements of the results overviewed in Section \ref{polymer o}. In Section \ref{SIM}, we introduce stochastic interfaces and their connections with random matrices, fleshing out the description seen in Section \ref{SIM o}. Finally, in Section \ref{thermodynamics}, we study the large-$N$ asymptotics of the interfaces seen in Section \ref{SIM}, expanding on the discussion in Section \ref{thermodynamics o}.

The following sections of the paper, Sections \ref{polymer proofs} through \ref{stproof}, are dedicated to proving results stated in Sections \ref{polymer} through \ref{thermodynamics}, and at times providing further details for computations made in these earlier sections. In particular, Section \ref{polymer proofs} and Section \ref{SIM proofs} are dedicated to proving the results stated in Section \ref{polymer} and Section \ref{SIM} respectively. The proofs and further details surrounding the discussion in Section \ref{thermodynamics} are spread across two sections, the latter of which, Section \ref{stproof}, is dedicated to our derivation of the surface tension of the bead model.  \\

{\em Acknowledgements.}  Research supported by the European Research Council 
(Grant Number 669306).

\section{Non-intersecting polymers} \label{polymer}

\subsection{Series and parallel inequalities}  \label{polymer 1} 
Recall the series bound \eqref{series bound} and parellel bound \eqref{parallel bound} first introduced in Section \ref{polymer o}. We now take a moment to give proofs of these two inequalities.

First we consider the series bound \eqref{series bound}, which we prove using concatenation in series as follows. Suppose $x,y,z$ are points in $\mathbb{Z}^2$, and $\pi = \{ a_0,\ldots,a_p\}$ is a 1-path from $x$ to $y$, and $\gamma = \{ b_0,\ldots, b_q \}$ is a $1$-path from $y$ to $z$. Then we define the concatenated $1$-path $\pi \oplus \gamma := \{ c_0, \ldots, c_p, c_{p+1},\ldots,c_{p+q} \}$ from $x$ to $z$ by setting $c_i = a_i$ for $i \leq p$ and $c_i := b_{i-p}$ for $i > p$. It is plain to check that this is indeed a path from $x$ to $z$, and that 
\begin{align*}
F(\pi \oplus \gamma ) = F(\pi) + F( \gamma). 
\end{align*}
Now suppose $\vec{x}, \vec{y}$ and $\vec{z}$ are $k$-points satisfying $\vx \leq \vy \leq \vz$. If $\pi$ is a $k$-path from $\vx$ to $\vy$ and $\gamma$ is a $k$-path from $\vy$ to $\vz$, then we may form a $k$-tuple $\pi \oplus \gamma := ( \pi_1 \oplus \gamma_1, \ldots, \pi_k \oplus \gamma_k )$ of paths from $\vx$ to $\vz$, which may or may not be intersecting. 
Recall that we say a $k$-point $(x_1,\ldots,x_k)$ is nice if either $x_{i+1}$ lies both strictly north and strictly west of $x_i$, or $x_{i+1} = x_i + \mathbf{e}_2$. We now prove the series bound \eqref{series bound}, which hinges on the observation that whenever the intermediate point $\vy$ is nice, then this resulting $k$-tuple $(\pi \oplus \gamma)$ is guaranteed to be non-intersecting.

\bp[Proof of the series bound \eqref{series bound}]

First we show that if $\vy$ is nice, then $\pi \oplus \gamma$ is a $k$-path---that is the sets $(\pi_i \oplus \gamma_i)$ and $(\pi_j \oplus \gamma_j)$ are disjoint for each $i \neq j$. Since each of $\pi$ and $\gamma$ are themselves $k$-paths, it is sufficient to show that for each $i \neq j$, the sets $\pi_i$ and $\gamma_j$ are disjoint. 

To see this, first consider the case $i < j$. In this case $y_j$ lies strictly north of $y_i$, and hence every point of $\gamma_j$ lies strictly to the north of every point of $\pi_i$, and hence $\pi_i \cap \gamma_j = \varnothing$. Alternatively, consider the case $i > j$. Since $y$ is nice, either $y_i$ lies both strictly north and strictly west of $y_j$, or $y_i$ is of the form $y_i = y_j + m \mathbf{e}_2$ for some positive integer $m$. In the former case, every point of $\gamma_j$ lies strictly east of every point of $\pi_i$, and hence $\pi_i \cap \gamma_j = \varnothing$. In the latter case, we must have $y_{j+1} = y_j + 1 \mathbf{e}_2$, and in particular, $\gamma_j$ must go through $y_j$ and then $y_j + 1 \mathbf{e}_1$ in order avoid $y_{j+1} = y_j + 1 \mathbf{e}_2$. It follows that $\gamma_j - \{ y_j\}$ consists of points that lie strictly east of $\pi_i$, and hence again $\pi_i \cap \gamma_j = \varnothing$. 

We have proved that all cases $\pi_i$ and $\gamma_j$ are disjoint, and hence the $k$-tuple $\pi \oplus \gamma$ of paths from $\vx$ to $\vz$ is a $k$-path. In particular, $\oplus$ is an injection from $\Gamma_{\vx \to \vy}$ to $\Gamma_{\vy \to \vz} \to \Gamma_{ \vx \to \vz  }$. Observing that the $k$-paths $\pi$ and $\gamma$ only overlap at $\vy$, we see that this operation is additive in the sense that
\begin{align} \label{fe}
F(\pi \oplus \gamma ) = F(\pi) + F( \gamma). 
\end{align}
Now suppose we have three $k$-points $\vx \leq \vy \leq \vz$ such that $\vy$ is nice. 
Using the positivity of the summands in \eqref{PF} with the fact that $\oplus$ is an injection to obtain the second line below, and \eqref{fe} to obtain the third, we have
\begin{align*}
Z_{ \vx \to \vz}(\beta) &:= \sum_{ \theta \in \Gamma_{ \vx \to \vz} } \exp \left( \beta F(\theta) \right)\\
&\geq \sum_{ (\pi,\gamma) \in \Gamma_{ \vx \to \vy} \times \Gamma_{ \vy \to \vz}  } \exp \left( \beta F( \pi \oplus \gamma) \right)\\
&= \sum_{\pi \in \Gamma_{\vx \to \vy } }   \exp \left( \beta F( \pi ) \right) \sum_{\gamma \in\Gamma_{ \vy \to \vz} }   \exp \left( \beta F( \gamma ) \right) \\
&=: Z_{\vx \to \vy} (\beta) Z_{ \vy \to \vz} (\beta),
\end{align*} 
establishing the series bound \eqref{series bound}.\ep

We now consider the dual procedure of concatenation in parallel, giving a proof of the parallel bound \eqref{parallel bound}.
 \bp[Proof of the parallel bound \eqref{parallel bound}]
Suppose we have $(k+j)$-points $\vx = (x_1,\ldots,x_{k+j})$ and $\vy = (y_1,\ldots,y_{k+j})$. Then every non-intersecting $(k+j)$-path $\pi$ gives rise to a unique decomposition
\begin{align*}
\pi = \pi' \boxplus \pi'',
\end{align*}
where $\pi'$ is a non-intersecting $k$-path from $\vx' = (x_1,\ldots,x_k)$ to $\vy' = (y_1,\ldots,y_k)$, and $\pi''$ is a non-intersecting $j$-path from $\vx'' = (x_{k+1},\ldots,x_{k+j})$ to $\vy'' = (y_{k+1},\ldots,y_{k+j})$. This procedure gives us an injection from $\Gamma_{\vx \to \vy}$ to $\Gamma_{ \vx' \to \vy' } \times \Gamma_{ \vx'' \to \vy''}$, and this injection is additive in the sense that
\begin{align} \label{FI}
F( \pi' \boxplus \pi'') = F(\pi') + F( \pi''). 
\end{align}
Using the positivity of the summands in \eqref{PF} with the fact that the map $\pi \mapsto \pi' \boxplus \pi''$ is an injection satisfying \eqref{FI}, we have
\begin{align*}
Z_{ \vx \to \vy}(\beta) &:= \sum_{ \pi \in \Gamma_{ \vx \to \vy }} \exp \left( \beta F(\pi) \right)\\
&\leq \sum_{\pi' \in \Gamma_{\vx' \to \vy' } }   \exp \left( \beta F( \pi' ) \right) \sum_{\pi'' \in\Gamma_{ \vx'' \to \vy''} }   \exp \left( \beta F( \pi'' ) \right) \\
&=: Z_{\vx \to \vy} (\beta) Z_{ \vy \to \vz} (\beta),
\end{align*} 
proving \eqref{parallel bound}.
\ep

We remark that by iterating \eqref{parallel bound}, for any pair of $k$-points $\vx = (x_1,\ldots,x_k)$ and $\vy = (y_1,\ldots,y_k)$ we obtain
\begin{align} \label{seperating}
Z_{\vx \to \vy}(\beta) \leq \prod_{i=1}^k Z_{x_i \to y_i}(\beta).
\end{align}
In fact, thanks to the celebrated Lindstr{\"o}m-Gessel-Viennot lemma \cite{sta}, the $k$-path partition functions have a determinantal expression in terms of the one-point partition functions:
\begin{align} \label{gv}
Z_{\vx \to \vy}(\beta) = \det_{i,j=1}^k \left( Z_{x_i \to y_j}(\beta) \right).
\end{align}
In light of \eqref{gv}, the inequality \eqref{seperating} may be understood as an analogue of Hadamard's inequality for the random matrix $(Z_{x_i \to y_j}(\beta))_{i,j=1}^k$.

\subsection{The single path partition function}
The case $k=1$ corresponds to studying the asymptotics of the partition function $Z_{(1,1) \to (n,m)}(\beta)$ as $m$ and $n$ become large, and is widely discussed in the literature. By a  standard subadditivity argument using the special case $k=1$ of the series bound \eqref{series bound} and Kingman's ergodic theorem \cite{kin}, it can be shown that the almost-sure limit
\begin{align} \label{infinite volume free energy}
f_c(\beta) := \lim_{N \uparrow \infty} \frac{1}{N} \log Z_{(1,1) \to (N,cN)}(\beta)
\end{align}
exists and is equal to $\sup_{N \geq 1} \frac{1}{N} \E \left[ \log Z_{(1,1) \to (N,cN)}(\beta)\right]$. 
Though there are useful bounds for $f_c(\beta)$ (see for instance Comets \cite{com}), explicit expressions for the its value remaining unknown but for a few cases which we now discuss. 

First we consider the infinite temperature limit case $\beta = 0$. In this case the partition function $Z_{(1,1) \to (n,m)}(0)$ is deterministic for every $(n,m)$, simply counting the number of paths starting at $(1,1)$ and ending at $(n,m)$. A straightforward computation using Stirling's formula tells us that the free energy $f_c(0)$ is given by
\begin{align*}
f_c(0) = \lim_{N \to \infty} \frac{1}{N} \log \binom{(c+1)N}{N}= (c+1)\log(c+1) - c \log c.
\end{align*}
We are also able to make sense of the zero temperature limit---the asymptotic case where $\beta \to \infty$. Here, the partition function concentrates on the path maximising the energy $F(\pi)$. Namely, for any pair of points $x \leq y$ in $\mathbb{Z}^2$, we have
\begin{align*}
\lim_{ \beta \uparrow \infty} \frac{1}{\beta} \log Z_{x \to y}(\beta) = \max_{\pi \in \Gamma_{x \to y}} F(\pi).
\end{align*}
In particular, without too much concern at this stage for the technical details surrounding the interchange of limits, we have 
\begin{align} \label{zero lim}
\lim_{\beta \uparrow \infty} \frac{1}{\beta} f_c(\beta) = \lim_{N \uparrow \infty} \frac{1}{N} \max_{\pi \in \Gamma_{(1,1) \to (N,cN)} } F(\pi) =: \ell_c
\end{align}
The value of $\ell_c$ is known in a few special cases. For example, in the case that each $\omega(z)$ is exponentially distributed with mean $1$ Rost \cite{ros} showed that 
\begin{align} \label{rost eq}
\ell_c = (1 + \sqrt{c})^2.
\end{align}
See \cite{oco:dir per} for further discussion of the single-path zero-temperature limit $\ell_c$.

Finally, there is one particular distribution for the random environment for which an explicit expression is known for $f_c(\beta)$ at a positive and finite value of $\beta$. Namely, in the case where each $e^{\beta \omega(z)}$ has the inverse-gamma distribution with parameter $\mu$ as in \eqref{lg1}, Sepp{\"a}l{\"a}inen \cite{sep} discovered a remarkable underlying algebraic structure based around the beta-gamma algebra making the partition function exactly solvable. (We refer the reader to \cite[Chapter 7]{com} for an overview.) Indeed, according to \cite[Theorem 2.4]{sep},
\begin{align} \label{sepp eq}
f_c(\beta) = - \sup_{\theta \in [0,\mu] } \left(  c \psi_0(\theta) + \psi_0( \mu  - \theta) \right),
\end{align}
where $\psi_0(\theta) := \Gamma'(\theta)/\Gamma(\theta)$ is the digamma function. Let us also mention here the papers \cite{OY,MO,CSS,OO,BC} where similar formulae have been obtained for some other exactly solvable polymer models.

\subsection{Asymptotics for finitely many long paths} \label{sec finite k}
For nice $k$-points $\vx$ and $\vy$, we study the asymptotics of the partition function $ Z_{\vx \to \vy + (N,cN)}(\beta)$ associated with the $k$-points stretched far apart in the $(1,c)$ asymptotic choice of direction. Here we are able to exploit the subadditivity due to the series bound \eqref{series bound} to prove the following results.

\begin{thm} \label{finite k thm}
There exists a function $f_c(k,\beta)$ such that for any pair of nice $k$-points $\vx$ and $\vy$, we have the almost sure convergence
\begin{align*}
\frac{1}{N} \log Z_{ \vx \to \vy + (N,cN)}(\beta) \to f_c(k,\beta).
\end{align*}
The limit satisfies $f_c(k,\beta) = \sup_{N \geq 1}\E \left[ \frac{1}{N} \log Z_{ \vx \to \vy + (N,cN)}(\beta)  \right]  = \lim_{N \to \infty}  \frac{1}{N}  \E\left[ \log Z_{\vx \to \vy + (N,cN) } \right]$, and this quantity is independent of the choices $\vx$ and $\vy$.
\end{thm}

We point out that by \eqref{infinite volume free energy} and Theorem \ref{finite k thm}, by definition we have $f_c(1,\beta) = f_c(\beta)$. We also remark that by the parallel bound \eqref{parallel bound}, it can be seen that the limits $f_c(k,\beta)$ satisfy the inequalities
\begin{align*}
f_c(k + j , \beta) \leq f_c(k ,\beta) + f_c(j ,\beta),
\end{align*}
for positive integers $k$ and $j$. 
We can say something much stronger however. According to the following theorem, $f_c(k,\beta)$ grows linearly in $k$. 

\begin{thm} \label{finite k thm 2}
The limits satisfy 
\begin{align*}
f_c(k,\beta) = k f_c(\beta).
\end{align*}
\end{thm}
 Finally, we will also prove the following result allowing comparison of $f_c(\beta)$ for different values of $c$.

\begin{thm} \label{ineq thm}
Let $0 < c < c'$ be positive reals. Then we have 
\begin{align*}
f_c(\beta) + (c' - c)  \beta \nu \leq f_{c'}(\beta) \leq \frac{c'}{c} f_c(\beta) - \left( \frac{c'}{c} - 1 \right) \beta \nu. 
\end{align*}
\end{thm}
Theorems \ref{finite k thm}, \ref{finite k thm 2} and \ref{ineq thm} are proved in Section \ref{polymer proofs}.


In the next section we look at the case where there are many paths of finite length.

\subsection{Asymptotics for many paths of finite length}
Where in the last section we had finitely many paths and let their lengths tend to infinity, in this section we do the opposite, considering many non-intersecting paths of fixed length running side-by-side, and letting the number of paths tend to infinity. In this direction, for $x \in \mathbb{Z}^2$ and $h = (h_1,h_2) \in \mathbb{Z}^2$, define the stacked $k$-point $\vx^{h \uparrow k}$ at $x$ in the $h$ direction by 
\begin{align} \label{dir stack}
\vx^{h \uparrow k}_i := x + (i-1)h, ~~ i = 1,\ldots,k.
\end{align}
For $x \leq y$ in $\mathbb{Z}^2$, and a direction $h$ such that $h_1 \leq 0 < h_2$, we now consider the the large $k$-asymptotics of the random variable $
Z_{\vx^{h \uparrow k} \to \vy^{h \uparrow k} }(\beta)$. 
The parallel bound \eqref{parallel bound} gives us a subadditivity in the $k$-variable which we use to prove the following result.
\begin{thm} \label{many k}
Let $x \leq y$ be points in $\mathbb{Z}^2$, and $h = (h_1, h_2)$ such that $h_1 \leq 0 < h_2$. Then as $k \to \infty$, the random variable $\frac{1}{k} \log Z_{\vec{x}^{h \uparrow k} \to \vy^{h \uparrow k} }(\beta)$ converges almost surely to a deterministic limit $I_{y -x,h}(\beta)$.  
\end{thm}
In fact, in the infinite temperature case ($\beta =0$), for certain choices of $h$, the asymptotic limit $I_{z,h}(0)$ can be computed explicitly by using the asymptotic theory of Toeplitz determinants \cite{BS}. In this direction we define the symbol associated with $z$ and $h$ to be the Laurent series $a_{z,h}: \mathbb{T} \to \mathbb{C}$ on the unit circle given by
\begin{align*}
a_{ z, h}( s ) := \sum_{ m \in \mathbb{Z}} Z_{0 \to z +m h }(0) s^m.
\end{align*}
Every continuous function $a: \mathbb{T} \to \mathbb{C}$ has a unique decomposition $a(e^{it}) = |a(e^{it})| e^{ i c(t)}$ where $c:[0,2\pi) \to \mathbb{R}$ is a continuous function satisfying $c(0) = 0$. The winding number of $a$ is the integer
\begin{align} \label{wind def}
\wind(a) := \frac{1}{ 2\pi} \lim_{ t \to 2 \pi} c(t).
\end{align}
Using Szeg{\"o}'s limit theorem for the asymptotics of Toeplitz determinants, in Section \ref{toe proof} we prove the following result, which gives an expression for $I_{z,h}(0)$ in terms of the symbol $a_{z,h}$. 
\begin{thm} \label{toeplitz thm} 
Suppose that $h_1 < 0 < h_2$ and the winding number of the symbol $a_{z,h}$ is zero. Then there is a Laurent series $\sum_{m \in \mathbb{Z}} c_m s^m$ satisfying
\begin{align*}
a_{z,h}(s) = \exp \left( \sum_{ m \in \mathbb{Z}} c_m s^m \right),
\end{align*}
such that
\begin{align*}
\lim_{ k \to \infty} \left(  e^{ -k  c_0} Z_{\vx^{h \uparrow k} \to \vy^{h \uparrow k} }(0) \right) = \exp \left( \sum_{m=1}^\infty m c_m c_{-m} \right). 
\end{align*}
In particular, $I_{z, h}(0) = c_0$. 
\end{thm} 
To see an example of this result in action, let $y= (3,2)$, $x = (0,0)$, and let $h = (-2,2)$. To first get a rough idea of the asymptotics using the iterated parallel bound \eqref{seperating}, we note that there are $\binom{3+2}{2} = 10$ paths from $x$ to $y$, from which it follows that $Z_{ 0^{ h \uparrow k} \to z^{ h \uparrow k} }(0) \leq 10^k$, and hence $I_{(3,2),(-2,2)}(0) \leq \log 10$.

We now obtain the exact asymptotics using Theorem \ref{toeplitz thm}. First we note that the symbol associated with $x,y$ and $h$ is given by 
\begin{align*}
a_{y-x,h} (s) := \sum_{ m \in \mathbb{Z} } \binom{ 5 }{ 3 + 2 m } s^m   = \frac{5}{s} + 10 + s.
\end{align*}
We remark that $\Realpart \left(a_{y-x,h}(e^{it}) \right) = 6 \cos (t ) + 10 $ for each $t$ in $[0,2 \pi)$. It follows that $a_{y-x,h}$ has positive real part on $\mathbb{T}$, and hence has winding number $0$. It is straightforward to show that\ the symbol has representation 
\begin{align*}
a_{y-x,h} (s) =\exp \left(  \log ( 5 + 2 \sqrt{5} ) + \sum_{m = 1}^{\infty } \frac{(-1)^m}{ m } \left(\frac{c}{10} s\right)^m +  \sum_{ m = 1}^{ \infty } \frac{(-1)^m }{ m } \left(\frac{c}{2} s^{-1} \right)^m \right).
\end{align*}
It then follows from Theorem \ref{toeplitz thm} that
\begin{align*}
I_{(3,2), (-2,2)} = \log ( 5 + 2 \sqrt{5} ) \approx \log 9.472 \leq \log 10.
\end{align*}
\subsection{Asympotics for many long paths}
Finally, in this section we consider the asymptotics of the partition functions associated with many paths in a rectangle whose dimensions are growing to infinity together with the number of paths. 
For $k \geq 1$ and $x \in \mathbb{Z}^2$, we define the stacked $k$-point $\vx^{\uparrow k}$ above and below a point $x$ by 
\begin{align} \label{stacked}
\vx^{\uparrow k}_i := x + (i-1)\mathbf{e}_2, ~~~ \vx^{\downarrow k}_i := x + (k - 1 + i)\mathbf{e}_2, ~~ i = 1,\ldots,k.
\end{align}
For integers $m,n,k$, consider $k$-points $(1,1)^{\uparrow k}$ and $(n,m)^{\downarrow k}$ at the bottom-left and top-right corner of an $m \times n$ rectangle. Then there is at least one $k$-path from $(1,1)^{\uparrow k}$ to $(n,m)^{\downarrow k}$ if and only if $k \leq m \wedge n$. With this picture in mind, we consider the scaling regime $m = cN, n = N, k = \alpha N$ for $0 < \alpha \leq c \leq 1$, leading us to study the asymptotic growth of the random variable
\begin{align} \label{scaled part}
\frac{1}{N^2} \log Z_{(1,1)^{\uparrow \alpha N}  \to (N,cN)^{\downarrow \alpha N}  }(\beta).
\end{align}
To understand the $\frac{1}{N^2}$ scaling, it is useful first to consider the infinite-temperature limit $\beta = 0$, which amounts to a computation calculating the asymptotic size of the set $\Gamma_{ (1,1)^{\uparrow k}  \to (n,m)^{\downarrow k}}$ of $k$-paths on the $m \times n$ rectangle.
Indeed, in Section \ref{infinite temperature}, we exploit determinantal identities to prove the following result.
\begin{lem} \label{infinite growth}
As $N \to \infty$, $\frac{1}{N^2}\log Z_{(1,1)^{\uparrow \alpha N}  \to (N,cN)^{\downarrow \alpha  N}  }(0) \to w(c,\alpha )$, where the limit $w(c,\alpha)$ is given by 
\begin{align} \label{w lim}
w(c,\alpha ) = Q(1+c - \alpha ) + Q(1 - \alpha ) + Q(c - \alpha ) + Q (  \alpha ) - Q ( c) - Q(c+1 - 2 \alpha ),
\end{align}
and $Q(u) = \frac{u^2}{2} \log u$. 
\end{lem}
With the scaled-$k$ high temperature limit at hand, in Section \ref{jensen} we use Jensen's inquality to prove the following result relating the positive temperature partition function to the infinite-temperature limit.  For points $a = (a^1,a^2)$ and $b = (b^1,b^2)$ of $\mathbb{Z}^2$, define $|| a - b||_1 := |a^1 - b^1| + |a^2 - b^2|$. 
\begin{lem} \label{infinite bounds}
Let $\vx$ and $\vy$ be $k$-points such that $\vx \leq \vy$, and let $\varpi( \vx, \vy) := \sum_{i=1}^k ||y_i - x_i||_1$ be the number of environment weights in a $k$-path from $\vx \to \vy$ (not including the starting points). Suppose that $\nu := \E[ \omega(0,0)]$ and $G(\beta) := \E[ e^{\beta \omega(0,0) }]$ are finite. Then
\begin{align*}
 \beta \nu \leq \frac{1}{ \varpi( \vx, \vy)  } \E \ln \frac{ Z_{\vx \to \vy} (\beta)}{ Z_{\vx \to \vy} (0 )}  \leq  \log G(\beta).
\end{align*}
\end{lem}
Now define the upper and lower limits
\begin{align}
R^-(c,\alpha;\beta)&:= \liminf_{N \to \infty} \frac{1}{N^2} \E \log Z_{(1,1)^{\uparrow \alpha  N}  \to (N,cN)^{\downarrow \alpha N}  }(\beta) ,\label{LL} \\
R^+(c,\alpha;\beta)&:= \limsup_{N \to \infty} \frac{1}{N^2} \E\log Z_{(1,1)^{\uparrow \alpha  N}  \to (N,cN)^{\downarrow \alpha N}  }(\beta) . \label{UL}
\end{align}
The following result is an immediately corollary of Lemma \ref{infinite growth} and Lemma \ref{infinite bounds}, showing that these lower and upper limits may be sandwiched within terms involving the infinite temperature limit.
\begin{cor} \label{correct scaling}
The asymptotic $k$-path partition functions satisfy
\begin{align*}
\alpha (1 + c - \alpha ) \beta \kappa \leq R^-(c,\alpha;\beta)- w_c(t)  \leq R^+(c,\alpha;\beta)- w_c(\alpha ) \leq \alpha ( 1 + c - \alpha ) G( \beta). 
\end{align*}
\end{cor}
However, we conjecture the following stronger result.
\begin{conj}
The lower and upper limits $R^-_c(\alpha ,\beta)$ and $R^+_c(\alpha ,\beta)$ agree, and moreover in this case
\begin{align*}
\frac{d}{d\alpha } R^+(c,\alpha;\beta)\Big|_{\alpha = 0} = f_c(\beta).
\end{align*}
\end{conj}


\section{Stochastic interfaces} \label{SIM}

\subsection{Stochastic interfaces} \label{SI}
Let $\Lambda$ be a finite subset of $\mathbb{Z}^2$ and let $\Lambda^* := \left\{ \langle x, y \rangle  \in \Lambda^2 :  \text{$y - x$ is equal to $\mathbf{e}_1$ or $\mathbf{e}_2$} \right\}$ be the set of directed edges in $\Lambda$. Suppose now $V$ is a convex function on $\mathbb{R}$ and $(W_x : x \in \Lambda)$ are weight functions. A stochastic interface is a random function $\phi:\Lambda \to \mathbb{R}$ whose law is given by $Z^{-1} \exp( - H_\Lambda[ \phi] ) \prod_{ x \in \Lambda} d \phi(x)$, where $H_\Lambda$ is the Hamiltonian
\begin{align*}
H_\Lambda[\phi] := \sum_{ \langle x,y \rangle \in \Lambda^*} V \left( \phi(y) - \phi(x) \right) + \sum_{ x \in \Lambda} W_x (\phi(x)),
\end{align*}
and $Z := \int_{ \mathbb{R}^{\Lambda}} \exp( - H_\Lambda[ \phi] ) \prod_{ x \in \Lambda} d \phi(x) $ is the interface partition function. We observe that if $D$ is a subset of $\Lambda$, then the marginal law of $\phi$ on $D$ may be obtained by integrating out $\phi$ on $\Lambda - D$. Namely, for any function $\lambda:D \to \mathbb{R}$ defined on a subset $D$ of $\Lambda$,  \begin{align*}
P \left( \phi(y) \in d \lambda(y) , ~ y \in D \right) = \frac{ g_{D}(\lambda) }{ Z} \prod_{y \in D} d \lambda(y),
\end{align*}
where $g_D : \mathbb{R}^D \to \mathbb{R}$ is the energy integral
\begin{align*}
g_D( \lambda ) := \int_{ \mathbb{R}^{ \Lambda} } \exp( - H_\Lambda[\phi]) \prod_{y \in D} \delta_{\lambda(y)} ( d \phi(y) ) \prod_{ x \in \Lambda -D} d \phi(x). 
\end{align*}

We will be particularly interested in stochastic interfaces on triangular and square subsets of $\mathbb{Z}^2$. Considering triangular sets first, let $T_N := \{ (i,j) : 1 \leq i \leq j \leq N \}$ be the triangle in $\mathbb{Z}^2$, let $D_N := \{( i,i) \in T_N  \}$ be the set of diagonal entries and let $E_N := T_N - D_N$ be the non-diagonal entries. When $\Pi = D_N$, we call the energy integral a pattern integral and write $g := g_{D_N}$ for short. In other words, a pattern integral is simply a function $g^V:\mathbb{R}^N \to \mathbb{R}$ given by
\begin{align} \label{PI}
g^V(\lambda_1,\ldots,\lambda_N) := \int_{\mathbb{R}^{T_N}} \exp \left ( - H_{T_N}[\phi] \right) \prod_{x \in E_N } d \phi(x) \prod_{i=1}^N \delta_{u_i}( d \phi(i,i) ).
\end{align}
where $d \phi(x)$ is Lebesgue measure and $\delta_u$ is the Dirac mass at $u$. 

Though the majority of surrounding literature on stochastic interfaces is restricted to models with symmetric interaction potentials, we will be most interested in the potentials associated with the exponential and bead interaction models, neither of which are symmetric. These are given by 
\begin{align}
\mathsf{exp}(u) := e^u~~~\text{and}~~~ \mathsf{bead}(u) := \infty \ind_{u > 0 },
\end{align}
respectively, and will appear in stochastic interfaces relating to random polymers and random matrices which we introduce in the following two sections. 

For a moment let the weight functions $(W_x : x \in T_N)$ be zero and consider interfaces defined on the triangle $T_N$ with the bead interaction $V(u) = \mathsf{bead}(u)$. It is straightforward to see that the functions $\phi:T_N \to \mathbb{R}$ for which the Hamiltonian $H_{T_N}[\phi]$ is finite are precisely the Gelfand-Tsetlin patterns, namely those functions satisfying the inequalities
\begin{align*}
\phi({i,j-1}) \leq \phi({i,j})  \leq \phi({i+1,j}) ~~ \text{for every $1 \leq i< j \leq N$}.
\end{align*}
Given a vector $(\lambda_1,\ldots,\lambda_N)$, we write $GT_N(\lambda)$ for the set of Gelfand-Tsetlin patterns $\phi:T_N \to \mathbb{R}$ satisfying $ (\phi({1,1}),\ldots,\phi({N,N})) = (\lambda_1, \ldots, \lambda_N)$. Clearly this set is empty unless $\lambda_1 \geq \ldots \geq \lambda_N$. In fact, it is well-known that the pattern integral $g^{\mathsf{bead}}$ associated with the $V(u) = \mathsf{bead}(u)$ is given by
\begin{align} \label{wdf}
g^{\mathsf{bead}}(\lambda) = \int_{\mathbb{R}^{T_N} } \prod_{ \langle x,y \rangle \in T_N^*} \ind_{ \phi(y) \leq \phi(x) } \prod_{ x \in E_N} d \phi(x)  \prod_{i = 1}^N \delta_{\lambda_i}( d \phi(i,i) ) = \Vol{ GT_N(\lambda) } = \frac{ \Delta_N(\lambda)}{ H(N)} ,
\end{align}
where
\begin{align} \label{van}
\Delta_N(\lambda) := \prod_{1 \leq i \leq j \leq N} (\lambda_i - \lambda_j) \ind_{ \lambda_1 \geq \ldots \geq \lambda_N },
\end{align}
 and $H(N) := \prod_{j=0}^{N-1} j!$ is the superfactorial. By using the representation \cite[Section I.3]{mac} of the product occuring in \eqref{van} as a determinant of a Vandermonde matrix:
\begin{align*}
\prod_{1 \leq i < j \leq N} (\lambda_i - \lambda_j) = \det_{i,j=1}^N (\lambda_i^{N-j} ),
\end{align*}
it is straightforward to prove the volume formula \eqref{wdf} by induction.

As for the exponential interaction potential $\mathsf{exp}$, the associated pattern integral 
\begin{align} \label{wf}
g^{\mathsf{exp}}(\lambda) := \int_{\mathbb{R}^{T_N}} \exp \left ( - \sum_{(x,y) \in T_N^* } e^{ \phi(y) - \phi(x) } \right) \prod_{x \in E_N } d \phi(x) \prod_{k=1}^N \delta_{\lambda_k}( d \phi({i,i}) )
\end{align}
is known as a Whittaker function (with parameter $0$).

We record the following lemma, which we will use shortly to connect the diagonal entries of stochastic interfaces with both the eigenvalue ensembles of random matrix theory as well as the so-called Whittaker measure.

\begin{lem} \label{integrate out}
Let $\phi^N:S_N \to \mathbb{R}$ be a stochastic interface with interaction potential $V$ such that the off-diagonal weight functions $(W_x : x \in E_N)$ are identically zero. Then the marginal law of the diagonal $(\lambda_1,\ldots,\lambda_N) := (\phi^N(1,1),\ldots,\phi^N(N,N))$ is given by
\begin{align*}
\frac{1}{Z} g^V(\lambda)^2 \exp \left( - \sum_{i = 1}^N W_{(i,i)}( \lambda_i) \right) \prod_{i=1}^N d \lambda_i.
\end{align*}
\end{lem}
\bp
Integrating out the off-diagonal variables $( \phi(i,j))_{1 \leq i < j \leq N }$ and $(\phi(i,j))_{1 \leq j < i \leq N }$ and using the definition \eqref{PI}, we obtain two powers of $g^V(\lambda)$. 
\ep
In the next section we discuss the Whittaker measure, a collection of random variables related to the partition functions of a random polymer which may be thought of in terms of a stochastic interface. 

\subsection{The Whittaker measure as a stochastic interface}
We will see now that for a random polymer with product weights distributed according to the inverse gamma law \eqref{lg1}, the partition functions have an expression in terms of a stochastic interface. In this direction, first recall the inverse gamma distribution with parameter $\mu$
\begin{align} \label{lg2}
I_\mu(ds) = \frac{ s^{ - \mu - 1} e^{ - 1/s} }{ \Gamma(\mu) } \ind_{s > 0} ds, 
\end{align}
and let $F_\mu(s) := \int_0^s I_\mu (du)$ be the associated distribution function. In order to consider the simultaneous behaviour of our partition functions as $\mu$ varies, we would like to take a coupling of our weight variables as functions of $\mu$. To this end, suppose $\left( U(z): z \in \mathbb{Z}^2 \right)$ are uniformly distributed random variables on the unit interval, and for $z \in \mathbb{Z}^2$ define the random variable $\zeta_\mu(z) := F^{-1}_\mu(U(z))$. Clearly $\zeta_\mu(z)$ is inverse gamma distributed with parameter $\mu$. 

Now for $1 \leq k \leq m \leq N$, define the path sets 
\begin{align*}
\Gamma^{N}(m,k) := \Gamma_{ (1,1)^{\uparrow k} \to (N,m)^{\downarrow k} }, ~~~~ \tilde{\Gamma}^N(m,k) := \Gamma_{ (1,1)^{\uparrow k} \to (m,N)^{\downarrow k} },
\end{align*}
and consider the associated partition functions 
\begin{align} \label{tau def}
\tau^{N}_\mu(m,k) := \sum_{ \pi \in \Gamma^{N}(m,k)  } \prod_{ z \in \pi} \zeta_\mu(z), ~~~~~ \tilde{\tau}^{N}_\mu(m,k) := \sum_{ \pi \in \tilde{ \Gamma}^{N}(m,k)  } \prod_{ z \in \pi} \zeta_\mu(z).
\end{align}
We emphasise that in contrast to the partition functions studied in Section 2, the initial points $\{ (1,i-1) : i = 1,\ldots,k \}$ of the $k$-paths $\pi$ in $\Gamma^{N}(m,k) $ or $\tilde{\Gamma}^N(m,k) $ are included in the weight products in \eqref{tau def}. Now consider the random function $\varphi^N_\mu:S_N \to \mathbb{R}$ defined by
\begin{align} \label{phi def}
\varphi^N_\mu({i,j}) =
\begin{cases}
\log \left( \frac{ \tau^N_\mu({N-j+i,i})}{ \tau^N_\mu({N-j+i,i-1}) }   \right) ~~ i \leq j.\\
\log \left( \frac{ \tilde{\tau}^N_\mu({N-j+i,i})}{ \tilde{\tau}^N_\mu({N-j+i,i-1}) }   \right) ~~ i \geq j. \\
\end{cases}
\end{align}
with the conventions that $\tau^N_\mu({m,0}) = \tilde{\tau}^N_\mu({m,0}) = 1$ for each $m$. By way of a diagram the reader may convince themselves that $
\tau^N_\mu({N,k}) = \tilde{\tau}^N_\mu(N,k)$ for every $k \leq N$, from which it follows that the overlapping definitions in \eqref{phi def} on the diagonal $i=j$ agree. We remark that \eqref{phi def} may be inverted to give
\begin{align} \label{phi inv}
\log \tau^N_\mu(m,k) = \sum_{i=1}^k \varphi^N_\mu(i,N-m +i).
\end{align}
The following result is a consequence of results in \cite{COSZ}, and is the key idea in this section, stating that $\varphi^N_\mu$ is a stochastic interface:
\begin{prop} \label{interface 2}
The random function $\varphi^N_\mu:S_N \to \mathbb{R}$ given by \eqref{phi def} is a stochastic interface on $S_N$ with interaction potential $V(x) = \mathsf{exp}(x)$ and weight functions $W_{(i,j)}(u) = \ind_{i = j} \mu u + \ind_{ i = j = N} e^{ - u}$. The interface partition function is given by $Z_\mu := \Gamma(\mu)^{N^2}$. 
\end{prop}
In Section \ref{interface proofs} we give an explanation for how Proposition \ref{interface 2} follows from \cite{COSZ}. We remark that using the definition $g^{\mathsf{exp}}$ in Lemma \ref{integrate out}, it follows from Theorem \ref{interface 2} that the marginal law of the diagonal $\left(\lambda_1,\ldots,\lambda_N \right) := \left( \varphi_\mu^N(1,1),\ldots,\varphi_\mu^N(N,N) \right)$ may be given in terms of the so-called Whittaker measure with constant parameter $\mu$:
\begin{align} \label{constant whittaker}
P^N_\mu( d \lambda) := \frac{1}{\Gamma(\mu)^{N^2}} \exp \left(  - e^{ - \lambda _N } - \mu \sum_{i =1}^N \lambda_i \right) g^{\mathsf{exp}}(\lambda)^2 \prod_{ i =1}^N d \lambda_i.
\end{align}

\subsection{Eigenvalue processes and $\mu \to 0$}
 For an $N \times N$ matrix $A := (A_{i,j})_{1 \leq i,j \leq N }$, let $A^{(k)} := (A_{i,j})_{1 \leq i ,j \leq k }$ be its principal $k \times k$ minor.  Now suppose $H$ and $U$ are $N \times N$ complex matrices, where $H$ is Hermitian and $U$ is unitary. We define the eigenvalue process $\varphi_{H,U}:S_N \to \mathbb{R}$  associated with $H$ and $U$ to be the function $\varphi_{H,U}:S_N \to \mathbb{R}$ given by
\begin{align} \label{kappa def}
\varphi_{H,U} (i,j) := 
\begin{cases}
\text{The $i^{\text{th}}$ largest eigenvalue of $H^{(N-j+i)}$}  ~~~~~~~~~~~~~\text{if} ~i \leq j,\\
\text{The $j^{\text{th}}$ largest eigenvalue of $(U^*HU)^{(N-i+j)}$} ~ ~~~\text{if}~ i \geq j.\\
\end{cases}
\end{align}
Let us remark that the $N \times N$ matrices $H$ and $U^* H U$ have the same eigenvalues, hence the two definitions in \eqref{kappa def} are consistent on the diagonal $i = j$. We also point out that thanks to Cauchy's interlacing theorem, the eigenvalues of each $A^{(k)}$ interlace those of $A^{(k+1)}$. This tells us that for any directed edge $\langle x, y \rangle$ in $S_N^*$, $\varphi_{H,U} (x) \geq \varphi_{H,U} (y)$, or in other words, that $\varphi_{H,U}$ is a Gelfand-Tsetlin pattern on the subset $T_N$ of $S_N$. \\ 

We now consider the eigenvalue processes associated with certain ensembles of unitarily invariant random matrices, introducing the Gaussian Unitary Ensemble (GUE) and the Laguerre Unitary Ensemble (LUE). We say a random $N \times N$ complex Hermitian matrix is GUE distributed if its law is given by
\begin{align} \label{GUE law}
P^N_{\mathsf{GUE}}(dH) := A_N \exp \left( - \frac{1}{2} \Tr (H^2) \right) dH,
\end{align}
and we say a random $M \times M$ complex Hermitian matrix is LUE distributed with underlying parameter $N$ if its law is given by 
\begin{align} \label{LUE law}
P^N_{\mathsf{LUE},M} (dH) := B_{N,M} \det(H)^{N-M} \exp \left( - \Tr(H) \right) \ind_{H \geq 0}dH,
\end{align}
where $\ind_{H \geq 0}$ is the indicator function on the set of positive definite matrices.
When $M = N$ in \eqref{LUE law}, we just say $H$ is LUE distributed for short. Since both the trace and determinant are invariant under unitary conjugations, it follows that if $H$ has either of the laws  \eqref{GUE law} or \eqref{LUE law}, then for any unitary $U$, $U^* H U$ has the same law as $H$.

We will also use the following property of LUE matrices, which is straightforward to prove using the representation of Laguerre matrices as a product of matrices with independent complex Gaussian entries \cite[Chapter 3]{for}.
\begin{rmk} \label{remark LUE} If $H$ is distributed with law $P^N_{\mathsf{LUE},N}$, then the marginal law of the minor $H^{(k)}$ is given by $P^N_{\mathsf{LUE},k}$. 
\end{rmk}
The following proposition, which states that the eigenvalue processes associated with GUE and LUE random matrices can be reformulated as stochastic interfaces, is a consequence of well known properties of the eigenvalues of unitarily invariant random matrices~\cite{bar,jn}.
\begin{prop} \label{eigen interface}
Suppose $H$ is GUE distributed and $U$ is Haar distributed on the set of unitary matrices. Then the eigenvalue process $\varphi_{H,U}$  associated with $(H,U)$ is a stochastic interface with interaction potential $\mathsf{bead}(u)$ and weight functions 
\begin{align*}
W^{\mathsf{GUE}}_{(i,j)}(u) := \ind_{i = j} \frac{1}{2} u^2. 
\end{align*}
Alternatively, suppose $H$ is LUE distributed and $U$ is Haar distributed on the set of unitary matrices. Then the eigenvalue process $\varphi_{H,U}$  associated with $(H,U)$ is a stochastic interface with interaction potential $\mathsf{bead}(u)$ and weight functions 
\begin{align*}
W^{\mathsf{LUE}}_{(i,j)}(u) := \ind_{i = j} u + \ind_{ i = j = N} \infty \ind_{ u < 0}. 
\end{align*}
\end{prop}

We provide a proof of Proposition \ref{eigen interface} in Section \ref{interface proofs}. For short, we refer to the eigenvalue processes occuring in Proposition \ref{eigen interface} as the GUE and LUE eigenvalue processes. With this result at hand (in fact for the moment we only require the part regarding LUE matrices), we are ready to establish the relationship between the partition functions of the polymer with inverse-gamma distributed weights and the eigenvalues of random matrices.

\begin{thm}  \label{small mu}
Let $\varphi_\mu^N:S_N \to \mathbb{R}$ be the interface defined in \eqref{phi def}. Then as $\mu \downarrow 0$, the rescaled process $\mu \varphi_\mu^N$ converges in distribution to the LUE eigenvalue process.
\end{thm}
\bp
By Theorem \ref{interface 2}, $\varphi_\mu^N$ is a stochastic interface with the exponential interaction potential and weight functions $W_{(i,j)}(u) = \ind_{ i = j} \mu u + \ind_{ i = j = N} e^{ - u}$. It is immediate that the change of variable $\mu \varphi^N_\mu$ is also a stochastic interface with interaction potential $V(u) = \mathsf{exp}( u/\mu)$ and weight function $W(i,j)(u) = \ind_{i=j} u + \ind_{i=j=N} \frac{ e^{-u/\mu}}{\mu}$. \\

Taking $\mu \downarrow 0$, we have the (Lebesgue-almost-everywhere) convergence of both the interaction potential $\mathsf{exp}( u/\mu)$ to $\mathsf{bead}(u)$ and the weight
$\frac{ e^{-u/\mu}}{\mu}$ to $\infty  \ind_{u < 0}$, which are precisely the interactions and weights associated with $\varphi_{\mathsf{LUE}}^N$.
\ep

We now consider the implications of this connection by looking directly at the small-$\mu$ asymptotics of the polymer. Suppose $\zeta_\mu$ is inverse-gamma distributed with parameter $\mu$, and let $\chi_\mu := \mu \log \zeta_\mu$. Using the small-$\mu$ asymptotics of the Gamma function, it is straightforward to verify that
\begin{align*}
\E[ e^{ - s \chi_\mu } ] = \frac{ \Gamma (\mu (s + 1 ) ) }{ \Gamma ( \mu) } \to \frac{1}{1+s} ~~~ \text{as $\mu \downarrow 0$}.
\end{align*}
That is, $\chi_\mu$ converges in distribution to a standard exponential random variable as $\mu \downarrow 0$. Applying these facts to the partition functions \eqref{tau def}, we see that as $\mu \downarrow 0$ we have the convergence in distribution
\begin{align} \label{go mu}
\mu \log \tau^N_\mu(m,k)=  \mu \log  \sum_{ \pi \in \Gamma^{N}(m,k)  } \exp \left( \frac{1}{\mu} \sum_{z \in \pi} \mu \log \zeta_\mu(z)\right) \to^{(d)} \max_{ \pi \in  \Gamma^{N}(m,k) } \sum_{z \in \pi} e(z) =: L^N(m,k),
\end{align}
where $( e(z) : z \in \mathbb{Z}^2)$ are independent standard exponential random variables.

On the other hand, as $\mu \downarrow 0$, combining \eqref{phi inv} and Theorem \ref{small mu}, we have the convergence in distribution
\begin{align} \label{go mu 2}
\mu \log \tau^N_\mu(m,k) := \sum_{i=1}^k \mu \varphi_\mu^N(i,N-m+i) \to \sum_{ i = 1}^k \varphi_{\mathsf{LUE}}^N(i,N-m+i).
\end{align}
Comparing \eqref{go mu} with \eqref{go mu 2} we obtain the distributional equality
\begin{align} \label{law eq}
L^N(m,k) =^{(d)} \sum_{i=1}^k \varphi_{\mathsf{LUE}}^N (i, N-m+i). 
\end{align}
By using the definition \eqref{kappa def} in \eqref{law eq} and considering Remark \eqref{remark LUE}, we recover the following result by Johansson \cite[Proposition 1.4]{joh} connecting the eigenvalues of LUE matrices with non-intersecting last passage percolation on an exponential polymer. (See also Doumerc \cite{dou}.)

\begin{prop}
The random variable $L^N(m,k)$ is equal in law to the sum of the $k$ largest eigenvalues of an $m \times m$ Laguerre matrix with underlying parameter $N$. 
\end{prop}

\subsection{The $\mu \to \infty$ limit} \label{inf finite}

Where in the last section we considered the small-$\mu$ asymptotics of the interface $\varphi_\mu^N$, here we consider the large-$\mu$ asymptotics. The following theorem states that as $\mu \to \infty$, a suitable rescaling of the $\varphi_\mu^N$ converges to a explicit deterministic shape related to a combinatorial problem. 

\begin{thm} \label{large mu}
Let $\theta_\mu^N:S_N \to \mathbb{R}$ be the rescaled interface 
\begin{align} \label{theta def}
\theta^{N}_{ \mu}(i,j):= \varphi^N_\mu(i,j) + (  2N  + 1 - (j + i))\log \mu .
\end{align}
Then as $\mu \to \infty$, $\theta^N_\mu$ converges almost-surely to a deterministic function $\theta_{\mathsf{min}}^N:S_N \to \mathbb{R}$, where $\theta_{\mathsf{min}}$ is the symmetric function given by
\begin{align} \label{psi form}
\theta^N_{\mathsf{min}}(i,j) := \log \left( \frac{ (i-1)! (2N - j - i + 1)! (2N - j - i)! }{(2N-j)! (N-j)! (N-i)! } \right), ~~~ i \leq j.
\end{align}
Moreover, the function $\theta_{\mathsf{min}}^N$ is the minimiser of the energy functional $\mathcal{F}: \mathbb{R}^{S_N} \to \mathbb{R}$ given by 
\begin{align} \label{cal f}
\mathcal{F}^N[ \theta ] :=  e^{ - \theta (N,N)} + \sum_{i=1}^N \theta(i,i) + \sum_{ \langle x, y \rangle \in S_N^* } e^{ \theta(y) - \theta(x) }.
\end{align}
\end{thm}
Theorem \ref{large mu} is proved in Section \ref{large mu proof} by directly analysing the large-$\mu$ behaviour of the interface $\varphi_\mu^N$, and studying the relationship  this interface has with a deterministic polymer via the application of the law of large numbers to $1/\zeta_\mu(z)$ (which for integer values of $\mu$ has a representation at the sum of $\mu$ independent exponential random variables). The fact that the function $\theta_{\mathsf{min}}^N$ minimises $\mathcal{F}^N[\theta]$ is an offshoot of our result, and we suspect there is a more direct combinatorial proof. 
That completes the section on finite interfaces. In the next section, we look at the asymptotics of these interfaces as $N \to \infty$.


\section{Scaling limits of stochastic interfaces} \label{thermodynamics}

\subsection{Surface tension and asymptotics of stochastic interfaces} \label{ST}
In this section we willl be interested in applying thermodynamic heuristics to study the macroscopic behaviour of the stochastic interfaces seen in the previous section, formulating their limiting shapes in terms of variational problems, and using these limit shapes to anticipate the asymptotics of the partition functions associated with the non-intersecting log-gamma polymer. 

We begin by introducing surface tension, an asymptotic measure of the cost of an interface to lie at a certain tilt. For further information on surface tension, the reader is referred to Sheffield \cite{she} and Funaki \cite{fun}. Let $\partial S_N$ and $S^{\circ}_N$ denote the interior and the boundary of the square $S_N$. Consider the square Hamiltonian
\begin{align*}
H_{S_N}[\phi] := \sum_{ \langle x,y \rangle \in S_N^* } V \left( \phi(y) - \phi(x) \right)
\end{align*}
associated with an interaction potential $V$ and $W_x \equiv 0$ for every $x$. For $N \geq 3$, we define the finite surface tension of interaction potential $V$ at tilt $(p,q) \in \mathbb{R}^2$ by 
\begin{align} \label{finite surface}
\sigma^V_N(p,q) := - \frac{1}{(N-2)^2}  \log \left(   \int_{\mathbb{R}^{\Lambda_N }   }  \exp ( - H_{S_N}[\phi] )   \prod_{ y \in \partial S_N } \delta_{(p y_1 + q y_2)}( d\phi(x) ) \prod_{ x \in S^{\circ}_N} d\phi(x)  \right).
\end{align}
The following result by Funaki and Spohn \cite{FS} states that under relatively strong conditions on the interaction potential $V$, the finite surface tensions $\sigma^V_N$ converge to a limit as a $N \to \infty$. The limit $\sigma^V$ is called the surface tension associated with $V$.
\begin{prop} \label{FS prop}
Suppose the potential function $V$ is symmetric, twice differentiable, and satisfies $
c_{-} \leq V''(x) \leq c_+
$
for some positive reals $c_{-} \leq c_+$. Then the limit $\sigma^V(p) := \lim_{N \to \infty} \sigma^V_N(p)$ 
exists and is a convex function on $\mathbb{R}^2$. 
\end{prop}

Neither of the interaction potentials $\mathsf{bead}$ and $\mathsf{exp}$ are symmetric, nor do they satisfy the technical condition required by Funaki and Spohn for the existence of the associated surface tension. Nonetheless, it was established by Sun \cite{sun} that we have the existence of the surface tension 
\begin{align*}
\sigma^{\mathsf{bead}}(p,q) := \lim_{N \to \infty} \sigma_N^{\mathsf{bead} }(p,q)
\end{align*}
associated with the $\mathsf{bead}(u)$ interaction, and moreover, Sun provided a formula for $\sigma^{\mathsf{bead}}(p)$ in a different coordinate system which we discuss below. Based on the relative similarity between the exponential and bead interaction, we anticipate the following conjecture.

\begin{conj} \label{conj st}
The surface tension $\sigma^{\mathsf{exp}} := \lim_{N \to \infty} \sigma^{\mathsf{exp}}_N$ associated with the exponential interaction exists.
\end{conj}
With the definition of surface tension at hand, we are ready to study the macroscopics of stochastic interfaces. First we define a rescaling of stochastic interfaces from $S_N$ to the unit square in $\mathbb{R}^2$.

\begin{defn} \label{rescaled}
Let $S = [0,1]^2$. The rescaled interface associated with a function $\phi^N:S_N \to \mathbb{R}$ is the function $\bar{\phi}^N : S \to \mathbb{R}$ given by 
\begin{align*}
\bar{\phi}^N(s,t) := \frac{1}{N} \phi ( s'N, t'N),
\end{align*}
where $s'$ is the smallest multiple of $1/N$ greater than $s$, and $t'$ is defined similarly.

\end{defn}

In Chapter 6 of \cite{fun}, Funaki shows the asymptotics of stochastic interfaces with potentials satisfying the conditions in Proposition \ref{FS prop} may be given in terms of minimisers of variational problems. Furthermore, it is known that the asymptotic shapes of interfaces with the bead interaction potential have intimate relationships with models in free probability  (see e.g.  Metcalfe \cite{met}). With these observations in mind, and in the event that Conjecture \ref{conj st} holds, we are lead to further predict the following result about the asymptotic shape of a certain class of interfaces.

\begin{conj} \label{main conj}
Let $\phi^N$ be a stochastic interface model with interaction potential either $V = \mathsf{bead}$ or $V = \mathsf{exp}$,  and weight functions of the form
\begin{align} \label{gen weight}
W_{(i,j)}^N(u) = \ind_{i=j} \mu u + \ind_{i=j=N} N f_N(u/N),
\end{align}
where $f_N$ be a function satisfying $f_N(u) \to \infty \ind_{ u \notin A}$ as $N \to \infty$, and $\mu$ is a positive real number. Suppose further that if $Z_N$ is the partition function associated with the stochastic interface on $S_N$, then the sequence $\{Z_N\}$ satisfies $\frac{1}{N^2} \log Z_N \to E \in \mathbb{R}$. 

Then the rescaled interface converges pointwise almost surely to a deterministic limit---that is, each $(s,t)$ in $S$, $\bar{\phi}^N(s,t)$ converges almost surely to a deterministic limit $\xi(s,t)$. Moreover, the limit function $\xi$ is the minimiser in $\mathcal{C}^1(S)$ of the functional $\mathcal{E}: \mathcal{C}^1(S) \to \mathbb{R}$ 
given by 
\begin{align*}
\mathcal{E}[ v ] := \int_S  \sigma^V \left( \nabla v(s,t) \right) ds dt + \mu \int_0^1   v(s,s) ds + \infty \ind_{ v(1,1) \notin A} + E.
\end{align*}
Finally, the value of the functional at the minimiser is given by $\mathcal{E}[\xi] = 0$. 
\end{conj}

\subsection{Asymptotics of the Whittaker measure}

Now we consider the implications of Conjecture \ref{main conj} for the function $\varphi_\mu^N$ defined in \eqref{phi def} and appearing as a stochastic interface in Proposition \ref{interface 2}, with a particular focus on what this tells us about large-$N$ asymptotics of the partition functions $\tau^N_\mu(m,k)$ under the scaling $m = cN$, $k = \alpha N$.

Indeed, we anticipate that the rescaled interface $\bar{\varphi}^N_\mu$ converges pointwise almost-surely to a deterministic limit $\xi_\mu$, and that  $\xi_\mu$ is the minimiser over all functions $v$ in $\mathcal{C}^1(S)$ satisfying $v(1,1) \geq 0$ of the energy functional $\mathcal{E}_\mu:\mathcal{C}^1(S) \to \mathbb{R}$ 
\begin{align} \label{main energy}
\mathcal{E}_\mu[v] := \int_S \sigma^{\mathsf{exp}}( \nabla v(s,t) ) ds dt + \mu \int_0^1 v(s,s) ds + \log \Gamma(\mu).
\end{align}

In this case, we are able to recover the large-$N$ asymptotics of the many-path partition functions associated with the inverse gamma polymer of parameter $\mu$. Namely, the large-$N$ asymptotics of \eqref{phi inv} suggest that 
\begin{align} \label{many path}
\lim_{N \to \infty} \frac{1}{N^2} \log \tau^N_\mu( cN, \alpha N) := \int_0^\alpha \xi_\mu(u, 1- c  +u) du. 
\end{align}

Short of offering an explicit expression for the limit shape $\xi_\mu$, we make a few predictions about its properties. First of all, by the symmetry of the energy functional $\mathcal{E}_\mu$ \eqref{main energy}, it is plain that the minimiser $\xi_\mu$ is itself symmetric. 

Moreover, consider minimising the two competing terms appearing in $\mathcal{E}_\mu[v]$ over functions satisfying $v(1,1) = 0$. On the one hand, the first term $\int_S \sigma^{\mathsf{exp}}( \nabla v)$ encourages $v$ to have negative derivatives with respect to both $s$ and $t$, where as the intermediate term $\mu \int_0^1 v(s,s) ds $ wants $v$ to decreases. As $\mu$ becomes larger, this second effect becomes stronger, and we anticipate that
\begin{align*}
\text{$\xi_\mu$ is monotone decreasing in $\mu$}. 
\end{align*}

With these broad observations made, we now turn to discussing existing results in the literature which give us the values of $\xi_\mu$ on the boundaries of the square $S$. First we note that by \eqref{phi def} the values taken by $\varphi^N_\mu$ on the boundary of $S_N$ are given in terms of the one-path logarithmic partition functions
\begin{align*}
\varphi_\mu^N(1,j) := \log \tau_\mu^N( N - j + 1, 1),
\end{align*}
which suggests 
\begin{align} \label{boundary tau}
\xi_\mu(0,t) = \lim_{N \to \infty} \frac{1}{N} \log \tau_\mu^N \left( N (1-t), 1 \right).
\end{align}
In particular, by Sepp{\"a}l{\"a}inen's equation \eqref{sepp eq} we have
\begin{align} \label{psi edge}
\xi_\mu(0, t) := - \sup_{ \theta \in [0, \mu] } \left(  (1-t) \psi_0(\theta) + \psi_0( \mu - \theta) \right),
\end{align}
where $\psi_0(\mu) := \Gamma'(\mu)/\Gamma(\mu)$ is the digamma function. Of course by the symmetry, $\xi_\mu(t,0)$ also satisfies \eqref{psi edge}.

We now seek to understand the values taken by $\xi_\mu$ on the other boundaries of the square $S$, namely points of the form $(t,1)$ and $(1,t)$. To this end, again by \eqref{phi def} we have
\begin{align*}
\varphi^N_\mu(i,N) := \log \left( \frac{ \tau_\mu^N ( i,i) }{ \tau_\mu^N(i,i-1) } \right).
\end{align*}
O'Connell and Ortmann study a random variable related to $\frac{ \tau_\mu^N ( i,i) }{ \tau_\mu^N(i,i-1) }$ in \cite{OO}. In particular, it follows from Equation (2.1) and Proposition 2.1 of \cite{OO} that we have the relation
\begin{align} \label{oo1}
\tau_\mu^N(i,i) / \tau^N_\mu( i , i -1) = 1/Z_\mu^{\mathsf{OO}}(n,m),
\end{align}
where for the choice $m = N-i+1$ and $n=i$, and $Z^{\mathsf{OO}}_\mu(n,m)$ is the partition function appearing in the introduction of \cite{OO}. Paraphrasing \cite[Theorem 1.1]{OO}, we have
\begin{align} \label{oo2}
\frac{1}{tN} \log Z^{\mathsf{OO}}_{N-tN+1,tN} \to \inf_{\theta > 0} \left( \frac{1}{t} \psi_0(\theta + \mu) - \psi_0(\theta) \right).
\end{align}
(We note here that this result was also obtained independently by Corwin, Sepp{\"a}l{\"a}inen and Shen \cite{CSS}.)
Combining \eqref{oo1} and \eqref{oo2} with the definition of $\xi_\mu(t,1) := \lim_{N \to \infty} \frac{1}{N} \varphi^N_\mu(tN,N)$ and using \eqref{phi def}, we obtain
\begin{align} \label{psi edge 2}
\xi_\mu(t,1) = \sup_{ \theta >0 } \left( t \psi_0(\theta) - \psi_0(\mu + \theta) \right). 
\end{align}
As a consistency check we remark that equations \eqref{psi edge} and \eqref{psi edge 2} agree that $\xi_\mu(0,1) = - \psi_0(\mu)$. 

Finally, we now show that the average of $\xi_\mu(s,t)$ along every diagonal $\{ (u,1-c +u) : u \in [0,c] \}$ is given by $-\psi_0(\mu)$. Namely, we show
\begin{align} \label{r3}
\frac{1}{c} \int_0^c \xi_\mu( u, 1-c  + u ) du = - \psi_0(\mu).
\end{align}
To see that \eqref{r3} holds, we note from \eqref{tau def} that $\tau^N_\mu(m,m) := \prod_{i = 1}^m \prod_{j = 1}^N \zeta_\mu({i,j})$ is simply a product of $m N$ independent inverse-gamma random variables with parameter $\mu$, and it follows from the law of large numbers that when $m = c N$, we have the almost sure convergence 
\begin{align} \label{r1}
\lim_{N \to \infty} \frac{1}{N^2} \log \tau_\mu^N(cN,cN)=  c \E[ \log \zeta_\mu({1,1}) ]. 
\end{align}
Now note that
\begin{align} \label{r2}
\E[ \log \zeta_\mu({1,1}) ] = \int_0^\infty \log s~ \frac{s^{- \mu - 1 } e^{ - 1/s} ds}{\Gamma(\mu)} = - \frac{\Gamma'(\mu)}{\Gamma(\mu) } = - \psi_0(\mu). 
\end{align}
Combining \eqref{r1} and \eqref{r2} with \eqref{many path}, we yield \eqref{r3}.

\subsection{The small-$\mu$ asymptotics of $\xi_\mu$ and the Mar\v{c}enko-Pastur law} \label{big N small mu}
We now consider the asymptotics of the limit shape $\xi_\mu:S \to \mathbb{R}$ as $\mu \to 0$. First we recall the Mar\v{c}enko-Pastur law concerning the asymptotic positions of eigenvalues of Laguerre matrices. \\

Suppose we have an $m \times m$ Laguerre distributed random matrix with underlying parameter $N$, and consider the asymptotics of the eigenvalues under the scaling limit $m =cN$ as $N \to \infty$. The Mar\v{c}enko-Pastur law states that the empirical measure of the $m$ eigenvalues of the rescaled matrix $\frac{1}{N} H$ converges almost surely to the Mar\v{c}enko-Pastur distribution $\nu_c$ with parameter $c$, where
\begin{align*}
\nu_c(du) := \frac{1}{2 \pi c} \frac{  \sqrt{  (M_c - u)(u - m_c) } }{ u } \ind_{u \in [m_c,M_c]} du,
\end{align*}
and $m_c = 1 + c - 2 \sqrt{c}$ and $M_c = 1 + c + 2 \sqrt{c}$. (There is also a version of the Mar\v{c}enko-Pastur distribution for $c > 1$ which we don't consider.)\\

We restate this result in terms of a position function $\rho_{\mathsf{mp}}:S \to \mathbb{R}$ as follows. Let $m =\lf c N\rf$ and $ k = \lf \alpha N \rf$ with $0 < \alpha \leq c \leq 1$, and define the random variable $\rho^N(c,\alpha)$ to be the the $k^{\text{th}}$ largest eigenvalue of a random matrix $\frac{1}{N} H$, where $H$ has law $P^N_{\mathsf{LUE},m}$. Then as $N \to \infty$, $\rho^N(c,\alpha)$ converges almost surely $\rho(c,\alpha)$, the  solution of the equation
\begin{align} \label{mp inv}
\int_{\rho_{\mathsf{mp}}(c, \alpha)}^{M_c}   \frac{1}{2 \pi c} \frac{  \sqrt{  (M_c - u)(u - m_c) } }{ u } \ind_{u \in [m_c,M_c]} du = \frac{\alpha}{c}.
\end{align}

Taking $c = 1 - t + s$ and $\alpha = s$, and using the definition of the eigenvalue process $\varphi_{\mathsf{LUE}}^N$, we yield the following proposition .

\begin{prop} \label{mp}
Let $\bar{\varphi}^N_{\mathsf{LUE}}:S \to \mathbb{R}$ be the rescaled interface associated with the Laguerre eigenvalue process. Then for each $(s,t)$ in $S$, $\bar{\varphi}_{\mathsf{LUE}}^N(s,t)$ converges almost surely to $\xi_{\mathsf{mp}}(s,t)$, where $\xi_{\mathsf{mp}}:S \to \mathbb{R}$ is the symmetric function satisfying
\begin{align} \label{psi mp}
\xi_{\mathsf{mp}}(s,t) := \rho_{\mathsf{mp}} \left( {1 - t + s}  , s \right) \text{   for $s \leq t$.}  
\end{align}
\end{prop}

On the other hand, assuming Conjecture \ref{main conj}, the representation of the eigenvalue process $\varphi_{\mathsf{LUE}}^N$ as an interface in \eqref{eigen interface} implies that the asymptotic limit $\xi_{\mathsf{mp}}$ is the minimiser of the energy functional
\begin{align*}
\mathcal{E}_{\mathsf{mp}}[v] := \int_S \sigma^{\mathsf{bead}}( \nabla v(s,t) ) ds dt + \int_0^1 v(s,s) ds + \infty \ind_{v(1,1) < 0 }.
\end{align*}
Reiterating the connection with the interface $\varphi_\mu^N$, Theorem \ref{mp} can be read as saying
\begin{align} \label{big con}
\xi_{\mathsf{mp}} = \lim_{N \to \infty} \bar{\varphi}_{\mathsf{LUE}}^N = \lim_{N \to \infty} \lim_{ \mu \to 0} \varphi_\mu^N,
\end{align}
where the convergence is pointwise almost sure on $S$. Assuming we can interchange the order of taking the limits $\mu \to 0$ and $N \to \infty$ in the final term of \eqref{big con}, we anticipate that
\begin{align} \label{small mu xi}
\mu \xi_\mu \to \xi_{\mathsf{mp}} ~ \text{as $\mu \to 0$}.
\end{align}

Finally, we now show that Proposition \ref{mp} implies Rost's equation \eqref{rost eq} for the asymptotics associated with last passage percolation on an exponential polymer. Namely, taking $\mu \downarrow 0$ in \eqref{boundary tau}, using \eqref{small mu xi} and \eqref{go mu} with $k=1$, we obtain
\begin{align*}
\xi_{\mathsf{mp}}(0,t) := \ell_{1-t},
\end{align*}
where $\ell_c$ is defined as in \eqref{zero lim}. It remains to note from \eqref{psi mp} and \eqref{mp inv} that 
\begin{align*}
\ell_c = \xi_{\mathsf{mp}}(0,1-c) = \rho_{\mathsf{mp}}( c , 0) = M_c = (1+ \sqrt{c})^2,
\end{align*}
as required. 

\subsection{The large-$\mu$ asymptotics of $\xi_\mu$} \label{big N big mu}

We now consider the implications of our discussion in Section $\ref{inf finite}$ as $N \to \infty$. Let $\theta_\mathsf{min}^N:S_N \to \mathbb{R}$ be the function defined in \eqref{psi form}, and let $\bar{\theta}_{\mathsf{min}}^N:S \to \mathbb{R}$ be the associated rescaling. By using Stirling's formula with the definition \eqref{psi form} of $\theta^N_{\mathsf{min}}$, it is possible to prove the following result. We provide details in Section \ref{big N big mu proof}.

\begin{thm} \label{big N big mu thm}
Let $q(u) = u \log u$. The limit $
\xi_{\mathsf{ht}}( s,t) := \lim_{ N \to \infty} \bar{\theta}_{\mathsf{min}}^N(s,t)
$
exists, and is given by the symmetric function satisfying 
\begin{align*}
\xi_{\mathsf{ht}}(s,t) = q(s) + 2q(2-s-t) - q(2-t) - q(1-t) - q(1-s) \text{   for $s \leq t$.}  
\end{align*}
\end{thm}
Again, assuming we can interchange the limits $N \to \infty$ and $\mu \to \infty$, by \eqref{theta def} we expect that 
\begin{align} \label{psi conv}
\lim_{ \mu \to \infty} \left( \xi_\mu (s,t)  + (2 - s - t) \log \mu \right) = \xi_{\mathsf{ht}}(s,t).  
\end{align}
Finally, in light of the fact that $\theta_{\mathsf{min}}^N$ minimises the energy functional \eqref{psi form}, it is natural to expect the rescaled limit interface  $ \xi_{\mathsf{ht}} :S \to \mathbb{R}$ is the minimiser of the energy functional
\begin{align*}
\mathcal{E}_{\mathsf{ht}}[v] := \int_S  \sigma^{\mathsf{ht}}( \nabla v(s,t) )   ds dt + \int_0^1 v(s,s) ds + \infty \ind_{ v(1,1) < 0 }, 
\end{align*}
where  $\sigma^{\mathsf{ht}}( \nabla v(s,t) )  = \exp \left( \frac{ \partial v}{ \partial s} \right) + \exp \left( \frac{ \partial v}{ \partial t} \right) $. We regard the energy functional $\mathcal{E}_{\mathsf{ht}}$ as an asymptotic analogue of the discrete functional in \eqref{cal f}. \\

\vspace{8mm}
We take a moment to recapitulate on the large and small-$\mu$ limits of $\xi_\mu$ we have seen, and on their related variational problems. The function $\xi_\mu$, defined as the limit as $N \to \infty$ of the rescaled interface associated with $\varphi_\mu^N$, is the minimiser of the energy functional
\begin{align*}
\int_S \sigma^{\mathsf{exp}}( \nabla v(s,t) ) ds dt + \mu \int_0^1 v(s,s) ds + \infty \ind_{ v(1,1) < 0 }.
\end{align*}
Anticipating that taking $\mu$ to either $0$ or $\infty$ commutes with taking $N \to \infty$, we predicted the following. On the one hand, as $\mu \to 0$, we expect that $\mu \xi_\mu$ converges to $\xi_{\mathsf{mp}}$, which has an explicit expression \eqref{psi mp} in terms of the Marcenko-Pastur distribution, and we predict is the minimiser of 
\begin{align*}
\int_S \sigma^{\mathsf{bead}}( \nabla v(s,t) ) ds dt + \int_0^1 v(s,s) ds +  \infty \ind_{ v(1,1) < 0 }.
\end{align*}
On the other hand, as $\mu \to \infty$, we expect that $\xi_\mu (s,t)  + (2 - s - t) \log \mu$ converges to $\xi_{\mathsf{ht}}$, which has the explicit expression \eqref{psi form} above, and we predict is the minimiser of 
\begin{align*}
 \int_S \sigma^{\mathsf{ht}}( \nabla v(s,t) )  ds dt + \int_0^1 v(s,s) ds + \infty \ind_{ v(1,1) < 0 }.
\end{align*}

\subsection{Gaussian fluctuations at high temperature} \label{gauss}
In the last section we looked at taking $\mu \to \infty$ after $N \to \infty$, anticipating the convergence to a limit shape. In this section, we develop a finer picture, taking $\mu$ and $N$ to infinity together through the scaling limit $\mu = \kappa N^2$, and using the central limit theorem to characterise the fluctuations of $\tau^N_\mu(m,k)$ in this regime.

Let $Q^N_{m,k}$ be the uniform law on $\Gamma^N(m,k)$, and let $\pi$ be a random variable with law $Q^N_{m,k}$. We make several assumptions about the asymptotic density of scaled $k$-paths in the $N \times N$ square. First, for $(u,v) \in S$ and $\alpha \leq c \leq 1$, we expect the existence of the limit
\begin{align} \label{lim dens}
 \lim_{N \to \infty} Q^N_{cN,\alpha N} \left(  (uN,vN) \in \pi \right) =: q_{c,\alpha}(u,v).
\end{align}
Moreover, for any $n$ distinct points $(u_1,v_1),\ldots,(u_n,v_n)$ in $S$ we anticipate that we have the asymptotic decoupling
\begin{align} \label{dec}
\lim_{N \to \infty} Q^N_{cN,\alpha N} \left( \{(u_1,v_1),\ldots,(u_n,v_n)\} \in \pi \right) = \prod_{i=1}^n q_{c,\alpha}(u_i,v_i)  .
\end{align}
This decoupling equation \eqref{dec} can be seen to hold, at least in a weak sense, via the well-known bijection between 
non-intersecting paths and lozenge tilings of a hexagon, and in fact one can write down an 
explicit formula for the limiting function $q_{c,\alpha}(u,v)$ via a formula given for the limit
shape of the corresponding tiling model in \cite{CLP}. 
We also refer the reader to the work of Johansson~\cite{joh2}, where the non-intersecting
paths model is related to an extended Hahn process.

Now let $\kappa > 0$, and consider the sequence of random processes $H^N_\kappa := \{ H^N_\kappa(s,t) : 0 \leq s \leq t \leq 1 \}$ given by 
\begin{align*}
H_\kappa^N(s,t) := k(N+m-k) \left( \log \kappa + 2 \log N  \right) + \log \tau_{\kappa N^2}^N(m,k) - \log \# \Gamma^N(m,k),
\end{align*} 
where $m = \lf tN \rf$ and $k = \lf s N \rf$. In Section \ref{gauss proof}, we sketch ideas leading us to the following conjecture. 

\begin{conj} \label{conj gauss} The limit \eqref{lim dens} exists and the satifies the asymptotic decoupling \eqref{dec}. Furthermore, the random process $\{ H^N_\kappa(s,t) : 0 \leq s \leq t \leq 1 \}$ converges in distribution to a centred Gaussian process $\{ H_{\kappa}(s,t) \}$ with covariance 
\begin{align*}
\E \left[H_\kappa(s,t) H_\kappa(\hat{s},\hat{t}) \right] = \frac{1}{\kappa} \langle q_{t,s}, q_{\hat{t},\hat{s}} \rangle_{L^2(S)},
\end{align*}
where $\langle  f, g \rangle_{L^2(S)} := \int_S f(u,v) g(u,v) du dv$ is the $L^2$ inner product on the unit square $S$.
\end{conj}

It is immediate from \eqref{lim dens} that $q_{c,c}(u,v) = \ind_{ v < c}$ in $L^2(S)$, and in particular, for any $t, \hat{t}$ in $[0,1]$ we have  $\langle q_{t,t}, q_{\hat{t},\hat{t} } \rangle = t \wedge \tilde{t}$. This observation implies that
\begin{align} \label{brownian edge}
\{  \sqrt{\kappa} H_{\kappa}(t,t) : 0 \leq t \leq 1 \} \text{ is a Brownian motion,}
\end{align}
though we mention here that is is possible to deduce \eqref{brownian edge} more directly by applying Donsker's invariance principle to the collection of random variables $\{ \tau^N_{\kappa N^2} (cN,cN) : c \in [0,1] \}$. 

\subsection{The semicircle law and surface tension in the bead model} \label{stbead}
In Section \ref{stproof} we combine tools from random matrix theory with our variational approach to macroscopics of stochastic interfaces to obtain the following explicit expression for the surface tension associated with bead model
\begin{align} \label{st}
\sigma_{\mathsf{tilted}}^{\mathsf{bead}}( p, q) := 
\begin{cases}
 - \log \left( |p| \cos \left( \pi \frac{ q}{ |p| } \right) \right) & \text{if} ~ p < 0, |q| < |p| \\
\infty &\text{otherwise},\\
\end{cases}
\end{align}
where $\sigma_{\mathsf{tilted}}^{\mathsf{bead}}$ is the surface tension in the change of coordinates given by $\sigma^V_{\mathsf{tilted}}( p,q) := \sigma^V  \left(  \frac{1}{2} p - q, \frac{1}{2}p + q \right)$.

An equivalent formula to \eqref{st} was obtained in Sun \cite{sun}, where it is proved by viewing the bead model as a continuous version of the Cohn-Kenyon-Propp \cite{CKP} dimer model. Sun uses Kasteleyn theory \cite{ken} to express the partition functions of this dimer model in terms of a contour integral, and studies the asymptotics of these contour integrals to obtain the expression \cite[Definition 5.4]{sun}.

Our approach, which we now overview, uses simpler technology from variational analysis and random matrix theory. On the one hand, we have the following result concerning the asymptotics of $\varphi_{\mathsf{GUE}}^N$, which is simply a restatement of the semicircle law of classical random matrix theory \cite{AGZ}.

\begin{thm}
The rescaled interface $\bar{\varphi}_{\mathsf{GUE}}^N$ associated with the eigenvalue process of a GUE random matrix converges pointwise almost surely to a symmetric function $\xi_{ \mathsf{sc}}:S \to \mathbb{R}$ satisfying
\begin{align} \label{psi sc}
\xi_{ \mathsf{sc}}(s,t) := \sqrt{ 1 - t + s } \rho_{\mathsf{sc}} \left( \frac{ s}{ 1 -t +s }\right) \qquad\text{for} \qquad s \leq t,
\end{align} 
where $\rho_{\mathsf{sc}}(x) \in [-2,2]$ is defined implicitly through the equation
\begin{align*}
\int_{\rho_{\mathsf{sc}}(x)}^2  \frac{1}{2\pi} \sqrt{ 4 - u^2 } du = x.
\end{align*}
\end{thm}

On the other hand, assuming the truth of a result analogous to Conjecture \ref{main conj}---with quadratic weight functions $\ind_{i = j } \frac{u^2}{2}$ in place of the linear weight functions $\ind_{i=j} \mu u$ occuring in \eqref{gen weight}---we expect that the limit $\xi_{\mathsf{sc}}$ is the minimiser of the energy functional
\begin{align} \label{sc funct} 
\mathcal{E}_{\mathsf{sc}}[v] := \int_S \sigma^{\mathsf{bead}}( \nabla v(s,t) ) ds dt + \frac{1}{2} \int_0^1 v(s,s)^2 ds. 
\end{align}

With these observations at hand, the main idea of our derivation is straightforward: if the minimiser of \eqref{sc funct} has the form \eqref{psi sc}, then the tilted surface tension of the bead model must be given by \eqref{st}. 

To sketch out a few further steps here, we develop the following scaling limit of the formula \eqref{wdf} for the volume of the Gelfand-Tsetlin polytope. Namely, suppose we have a strictly decreasing $\rho:[0,1] \to \mathbb{R}$. Then taking $N \to \infty$ in the definition of the Vandermonde determinant $\Delta_N(\lambda) := \prod_{1 \leq i < j \leq N} (\lambda_i - \lambda_j)_+$, we obtain
\begin{align} \label{van asy}
\bm{\Delta}[ \rho] &:= \lim_{N \to \infty} \left( \frac{1}{N^2} \log \Delta_N( N \rho(1/N),\ldots, \rho(N/N) )  - \frac{1}{2} \log N \right)\\
&= \int_{0 < s < t < 1 } \log \left( \rho(s) - \rho(t) \right) ds dt.
\end{align}
Moreover, the large-$N$ asymptotics of the superfactorial are given by 
\begin{align} \label{bar asy}
\log H(N) = N^2 \left( \frac{1}{2} \log N - \frac{3}{4} + o(1) \right).
\end{align}
(See e.g. Chen \cite{che}.) Combining \eqref{van asy} and \eqref{bar asy}, we see that on the one hand, since $g^{\mathsf{bead}}(\lambda) = \Lambda_N(\lambda)/ H(N)$, 
\begin{align} \label{comb asy}
\lim_{N \to \infty} \frac{1}{N^2} g^{\mathsf{bead}}( N \rho(1/N),\ldots, N \rho(N,N) )  = \int_{0 < s < t < 1 }  \log \left( \rho(s) - \rho(t) \right) ds dt + \frac{3}{4} 
\end{align} 
On the other hand, reading directly off the integral representation of the pattern integral $g^{\mathsf{bead}}$ as it appears in \eqref{wdf}, it is natural to expect that as defined,
\begin{align} \label{direct hydro}
\lim_{N \to \infty} \frac{1}{N^2} g^{\mathsf{bead}}( N \rho(1/N),\ldots, N \rho(N,N) ) = -   \min_{ v \in \mathcal{C}^1(T)_\rho  } \int_{T } \sigma^{\mathsf{bead}} (\nabla v),
\end{align}
where $T := \{ 0 \leq s \leq t \leq 1 \}$, and $\mathcal{C}^1(T)_\rho $ is the set of $\mathcal{C}^1$ functions on $T$ satisfying $v(s,s) = \rho(s)$ for all $s \in [0,1]$. 

Combining \eqref{comb asy} with \eqref{direct hydro}, we obtain the following thermodynamic analogue of the Gelfand-Tsetlin volume formula 
\begin{align} \label{weyl hydro}
\mathcal{M}^V[\rho] := \min_{ v \in \mathcal{C}^1(T)_\rho  } \int_{T } \sigma^{\mathsf{bead}} (\nabla v) =  - \int_{0 < s < t < 1} \log ( \rho(s) - \rho(t) ) ds dt - \frac{3}{4}.
\end{align}

We consider the functional derivatives of the functional $\mathcal{M}^V[\rho]$, showing that if $v^*$ is the minimiser of the integral over $T$ in \eqref{weyl hydro}, then $v^*$ satisfies the equation
\begin{align} \label{cauchy}
\frac{ \partial \sigma}{ \partial \tau}  \left( \frac{ \partial v^*}{ \partial r}(r,0) , \frac{ \partial v^*}{ \partial \tau}(r,0) \right) = \int_0^1 \frac{ 1}{ \rho(r) - \rho(s) } ds .
\end{align}
However, we also know that $v^*$ must also be related to the semicircle law as it appears in \eqref{psi sc}. By plugging $\xi_{\mathsf{sc}}$ into the differential equation \eqref{cauchy} , and using a homogoneity property of the surface tension, we determine $\sigma_{\mathsf{tilted}}^{\mathsf{bead}}$ as having the form in \eqref{st}. The full details of this argument are given in Section \ref{stproof}, though it bares remarking here that in pinciple we could have alternatively used the Mar\v{c}enko-Pastur distribution in place of the semicircle law, and we only opt to use the latter because the calculations involved are more straightforward.


\section{The Finite $k$ case} \label{polymer proofs}
This section is dedicated to proving the results stated in Section \ref{polymer}. In Section \ref{jensen} we prove Lemma \ref{infinite bounds}, which relates the partition function $Z_{\vx \to \vy}(\beta)$ to its infinite temperature counterpart, $Z_{\vx \to \vy}(0)$. Sections \ref{5.1} and \ref{5.2} are dedicated to giving proofs of Theorems \ref{finite k thm} and \ref{finite k thm 2} concerning the asymptotics of the partition functions associated with a fixed number of long non-intersecting paths. In Section \ref{5.3} we prove Theorem \ref{many k} and in Section \ref{toe proof} we use Szeg{\"o}'s limit theorem to prove Theorem \ref{toeplitz thm}. Finally, in Section \ref{infinite temperature} we give a derivation of Macmahon's formula leading to a proof of Lemma \ref{infinite growth}.

We will use the following notation. If $(X_n)_{n \geq 1}$ are random variables and $\theta$ is a real number, we will write
\begin{align*}
X_n \to \theta \qquad \text{a.s.e } \qquad \text{as     $n \to \infty$}
\end{align*}
to denote that the sequence $X_n$ converges both almost-surely and in expectation to $\theta$ as $n \to \infty$. Namely,
\begin{align*}
\P \left( X_n \to \theta \right) = 1 \qquad \text{and} \qquad  \E \left[ |X_n - \theta | \right] \to 0 \qquad \text{as $n \to \infty$}.
\end{align*}

\subsection{The infinite temperature sandwich bounds} \label{jensen}
Recall that $\varpi( \vx, \vy) := \sum_{ i = 1}^k || y_i - x_i ||_1$, $\nu := \E \left[ \omega(0,0) \right]$, and $G(\beta) := \E \left[ e^{ \beta \omega(0,0)} \right]$.We now use  Jensen's inequality to prove Lemma \ref{infinite bounds}, which states that we have the bounds
\begin{align} \label{jensen sandwich}
 \beta \nu \leq \frac{1}{ \varpi( \vx, \vy)  } \E \ln \frac{ Z_{\vx \to \vy} (\beta)}{ Z_{\vx \to \vy} (0 )}  \leq  \log G(\beta).
\end{align}

\bp[Proof of Lemma \ref{infinite bounds}]
The infinite temperature partition function $Z_{\vx \to \vy} (0 )$ counts the set $\Gamma_{\vx \to \vy}$ of $k$-paths from $\vx \to \vy$. Let $Q_{\vx \to \vy}$ be the uniform measure on $\Gamma_{\vx \to \vy}$, and let $\Pi$ be a k-path-valued random variable with law $Q_{\vx \to \vy}$. By definition,
\begin{align*}
 \frac{ Z_{\vx \to \vy} (\beta)}{ Z_{\vx \to \vy} (0 )} = Q_{ \vx \to \vy} \left[ e^{\beta F(\Pi)} \right],
\end{align*}
where $F(\pi) := \sum_{i=1}^k \sum_{x \in \pi_i - \{ x_i \} } \omega(x)$. In particular,  
\begin{align*}
\frac{1}{ \varpi( \vx, \vy)  } \P \ln \frac{ Z_{\vx \to \vy} (\beta)}{ Z_{\vx \to \vy} (0 )}  = \frac{1}{ \varpi( \vx, \vy)  } \P \ln Q_{ \vx \to \vy} \left[ e^{\beta F(\Pi)} \right].
\end{align*}
Under the measure $\P Q_{\vx \to \vy}$, $F(\Pi)$ is equal in law to the sum of $\varpi( \vx, \vy)$ independent and identically distributed random variables with law $\omega(0,0)$. 

Now on the one hand, using Jensen's inequality to interchange the order of $\P$ and $\ln$, we obtain
\begin{align*}
\frac{1}{ \varpi( \vx, \vy)  } \P \ln Q_{ \vx \to \vy} \left[ e^{\beta F(\Pi)} \right] &\leq \frac{1}{ \varpi( \vx, \vy)  }  \ln  \P Q_{ \vx \to \vy} \left[ e^{\beta F(\Pi)} \right] \\
&= \frac{1}{ \varpi( \vx, \vy)  }  \ln \left( G(\beta)^{ \varpi( \vx, \vy)} \right),
\end{align*}
yielding the upper bound in \eqref{jensen sandwich}. On the other hand, using Jensen's inequality to interchange the order of $\ln$ and $Q_{\vx \to \vy}$ we have
\begin{align*}
\frac{1}{ \varpi( \vx, \vy)  } \P \ln Q_{ \vx \to \vy} \left[ e^{\beta F(\Pi)} \right] &\geq \frac{1}{ \varpi( \vx, \vy)  }   \P Q_{ \vx \to \vy} \left[ \ln  e^{\beta F(\Pi)} \right] \\
&= \frac{1}{ \varpi( \vx, \vy)  } \left( p( \vx, \vy) \nu \beta  \right),
\end{align*}
giving the lower bound in \eqref{jensen sandwich}.
\ep

In a later proof, we will require the following blunt lower bound on $\E[ \log Z_{ \vx \to \vy} ( \beta)]$, which is an immediate corollary of \eqref{jensen sandwich}. 

\begin{cor} \label{bluntbound}
Let $\vx \leq \vy$. Then
\begin{align*}
\beta \nu \varpi( \vx, \vy) \leq \E \left[ \log Z_{ \vx \to \vy} (\beta) \right]
\end{align*}
\end{cor}
\bp
Since $\vx \leq \vy$, $\log Z_{ \vx \to \vy}( 0) \geq 0$. Now rearrange the lower bound in \eqref{jensen sandwich}.
\ep

\subsection{Diagonal $k$-points} \label{sec:diagonal}
For technical reasons, in our proofs of asymptotic results we will like to use a particular type of nice $k$-point with several desirable properties. Namely, for $x \in \mathbb{Z}^2$, define the diagonal $k$-point at $x^{ \Uparrow k} = ( x^{ \Uparrow k}_1, \ldots, x^{ \Uparrow k}_k )$ with base $x$ by
\begin{align} \label{diagonal kpoints}
x^{ \Uparrow k}_i := x + (i-1)( - \mathbf{e}_1 + \mathbf{e}_2 ).
\end{align}
Recall that for $1$-points $x = x_1  \mathbf{e}_1 +  x_2 \mathbf{e}_2$ and $y = y_1 \mathbf{e}_1 + y_2 \mathbf{e}_2$ of $\mathbb{Z}^2$,  we say $x \leq y$ if $x_1 \leq y_1$ and $x_2 \leq y_2$.  We have $x \leq y$ if and only if there is a $1$-path from $a$ to $b$. The following proposition gives us a useful property of diagonal $k$-points.

\begin{prop} \label{diagonal ordering}
Let $x$ and $y$ be $1$-points of $\mathbb{Z}^2$. Then $x^{\Uk} \leq y^{\Uk}$ if and only if $x \leq y$. 
\end{prop}
\bp
On the one hand, if $x^{\Uk} \leq y^{\Uk}$, then plainly there is a $1$-path from $x = x^{\Uk}_1$ to $y = y^{\Uk}_1$. On the other hand, if $x \leq y$, then there exists a $1$-path $\pi_1 = \{ x , z_1,\ldots,z_{p-1}, y \}$ from $x$ to $y$. Now for $i = 2,\ldots,k$, define
\begin{align*}
\pi_i := \left\{ x + (i-1) ( - \mathbf{e}_1 + \mathbf{e}_2) , z_1 +  (i-1) ( - \mathbf{e}_1 + \mathbf{e}_2 ) , z_2 +  (i-1) ( - \mathbf{e}_1 + \mathbf{e}_2 ) , \ldots, y +  (i-1) ( - \mathbf{e}_1 + \mathbf{e}_2) \right\}.
\end{align*}
Then $\pi = (\pi_1,\ldots, \pi_k)$ is a $k$-path from $x^{\Uk}$ to $y^{\Uk}$.
\ep

\begin{prop} \label{diagonal bounding}
Let $\vx = (x_1,\ldots,x_k)$ be a $k$-point with strictly increasing vertical coordinates in the sense that writing $x_i = (a_i,b_i)$, we have $b_1 < \ldots < b_k$. (In particular, all nice $k$-points have this property.) Then there exists a $1$-point $y_*$ in $\mathbb{Z}^2$ such that $y^{\Uk} \geq \vx$ for all $y \geq y_*$. Likewise there exists a $1$-point $z_*$ such that $z^{\Uk} \leq \vx$ whenever $z \leq z_*$. 
\end{prop}
\bp
With $r := \max\{ a_i : 1 \leq i \leq k  \}$, define $a'_i := r + k - i$ and $b_i' = b_i$. Then there is a horizontal $k$-path travelling from $\vx$ to  $\vx' := (x'_1,\ldots,x'_k)$ defined by $x'_i:=(a'_i,b_i')$, and hence $\vx \leq \vx'$. Now define $s := \max \{ b_i : 1 \leq i \leq k  \}$. Then there is a vertical $k$-path travelling from $\vx'$ to $( r+k - 1, s - k + 1)^{\Uk}$. Setting  $y_* := ( r+k - 1, s - k + 1)$, Proposition \eqref{diagonal ordering} ensures that $\vx \leq y^{\Uk}$ for all $y \geq y_*$. The proof of the existence of $z_*$ is similar. 
\ep

\subsection{Proof of Theorem \ref{finite k thm}} \label{5.1}
In this section we use a subadditivity argument to prove Theorem \ref{finite k thm}, which we recall states that there exists a function $f_c(k,\beta)$ such that for any pair of nice $k$-points $\vx$ and $\vy$, we have the almost sure convergence $\frac{1}{N} \log Z_{ \vx \to \vy + (N,cN)}(\beta) \to f_c(k,\beta)$. 

The main idea here is that by the series bound \eqref{series bound}, the random variables $Y_{\vx \to \vy} := - \log  Z \left[ \vx \to \vy \right]$ are subadditive with respect to the partial ordering $\leq$ of nice points $\mathbb{Z}^{2 \times k}$ in the sense that for any three nice $k$-points $\vx \leq \vy \leq \vz$, we have
\begin{align} \label{sub}
Y_{\vx \to \vz}  \leq Y_{ \vx \to \vy}  + Y_{ \vy \to \vz} .
\end{align}
We want to exploit this subadditivity to study the asymptotics of the random variable $Z_{ \vx \to \vy + (N,cN) }$ as $N$ tends to infinity. To do so, we recall Kingman's subadditive ergodic theorem \cite{kin}, giving a slight restatement of a later version of the theorem by Liggett \cite{lig}. (In the statement in \cite{lig}, the random variables are indexed over pairs of natural numbers, where as we allow our random variables to be indexed by the integers. 
It is not difficult to see the consequences of the theorem continue to hold in this setting.)

\begin{thm}  \label{kin thm}
Let $(X_{r,s})_{r < s}$ be a family of random variables indexed over either $\{ (r,s) \in \mathbb{Z}^2 : r < s \}$ or $\{ (r,s) \in \mathbb{Z}_{\geq 0}^2 : r < s \}$ and satisfying the following three conditions.
\begin{enumerate}
\item For every triple of integers $r < s < t$, we have $X_{r,t} \leq X_{r,s} + X_{s,t}$. 
\item For every $a \geq 0$, the joint distribution of the processes $(X_{r+a,s+a})_{0 \leq r < s}$ are the same as those of $(X_{r,s})_{0 \leq r < s}$.
\item The expectation $g_r := \E[X_{0,r}]$ exists, and there exists a real number $L$ such that $\frac{g_r}{r} \geq L$ for every $r \geq 1$.
\end{enumerate}
Then there is a constant $\theta$ such that for every $r$ in $\mathbb{Z}$,
\begin{align*}
X_{r,s}/s \to \theta \qquad \text{a.s.e} \qquad \text{as $s \to \infty$}.
\end{align*}
Moreover, for any $r$, $\theta = \inf_{s > 0} \frac{1}{s} \E[ X_{0,s}]$. 
\end{thm}

With Kingman's subadditive ergodic theorem now stated, the proof of Theorem \ref{finite k thm} is split into three parts. First, in Lemma \ref{special points} we prove that the conclusion of Theorem \ref{finite k thm} holds for rational $c$ and for all diagonal $k$-points situated on a certain lattice. Thereafter we prove Lemma \ref{key equations}, which gives us a pair of inequalities for the asymptotics of partition functions in a more general setting. Using Lemma \ref{key equations}, we deduce that the dependence on the parameter $c$ is continuous, and use this observation to prove that the conclusion of Theorem \ref{finite k thm} holds for all nice $k$-points and all positive slopes $c$. 

\begin{lem} \label{special points}
Let $c > 0$ be a rational number and let $p$ and $q$ be positive integers such that $q/p = c$. For integers $r \in \mathbb{Z}$ define the $k$-points $x_r := (rp,rq)^{\Uparrow k}$, and for pairs of integers $r < s$, define the random variables 
\begin{align} \label{X def}
X_{r,s} := - \frac{1}{p} \log Z_{x_r \to x_s} ( \beta).
\end{align}
Then there is a real number $f_c(k,\beta)$ such that for any $r$,  the random variables $ \frac{1}{s} X_{r,s}$ converge to $- f_c(k,\beta)$ a.s.e. as $s \to \infty$. Moreover $f_c(k,\beta)$ is independent of the choice of positive integers $p$ and $q$ satisfying $q/p = c$.
\end{lem}

\bp

First we show that with $X_{r,s}$ defined as in \eqref{X def}, we are in the set up of Theorem \ref{kin thm}. First of all, the inequality $X_{r,t} \leq X_{r,s} + X_{s,t}$ is a consequence of \eqref{sub}. The second condition is immediate from the definition of $X_{r,s}$ (using the fact that the polymer weights are independent and identically distributed). Finally, the fact that the final condition is satisfied is a consequence of Lemma \ref{infinite growth}. Indeed, by using the upper bound in Lemma \ref{infinite growth} to obtain the inequality below, we have
\begin{align} \label{zero bounds 1}
\E[ X_{0,r} ] &= - \frac{1}{p} \E \left[ \log Z_{x_0 \to x_r} ( \beta)  \right] \nonumber \\
&\geq  - \frac 1 p \left \{ \varpi\left( x_0, x_r \right) \log G(\beta) + \E \left[   \log Z_{x_0 \to x_r} ( 0)  \right] \right \} \nonumber\\
&= - r k \frac{ p + q}{ p } \log G(\beta)  - \frac{1}{p} \E \left[   \log Z_{x_0 \to x_r} ( 0)  \right],  
\end{align}
where we used the fact that $\varpi\left( x_0, x_r \right) = rk(p+q)$ to obtain the final equality above. Now using the parallel bound \eqref{parallel bound} with $\beta=0$ to prove the first inequality below, and the inequality $ \binom{b}{a} \leq \left( \frac{eb}{a}\right)^a$ to obtain the second, we have
\begin{align} \label{zero bounds}
\# \Gamma_{(0,0)^{\Uk} \to (rp,rq)^{\Uk} } &\leq k \# \Gamma_{(0,0) \to (rp,rq) } = k \binom{ r(p+q)}{ rp} \nonumber \\
&\leq k \left( \frac{ e(p + q) }{ p} \right)^{rp} 
\end{align} 
Plugging \eqref{zero bounds} into \eqref{zero bounds 1} we obtain
\begin{align*}
\E [ X_{0,r} ] \geq  - \frac{1}{ p } \log k - r \left( - \log \left( \frac{p+q}{p} \right) - k \frac{p+q}{p} G(\beta) - 1\right),
\end{align*}
which proves the third and final condition holds in Theorem \ref{kin thm}. 

Consequently, by Theorem \ref{kin thm}, for each $r \in \mathbb{Z}$, the limit
\begin{align*}
f_{p,q}(k, \beta) := - \lim_{s \to \infty} \frac{1}{s} X_{r,s} =  \frac{1}{p}  \lim_{s \to \infty} \frac{1}{s} \log Z_{ x_r \to x_s}( \beta)
\end{align*}
exists almost-surely as a constant not depending on $r$.

It remains to prove $f_{p,q}(k,\beta)$ is independent of the choice of $p$ and $q$. To see this, first note that we may assume without loss of generality that $r = 0$ since $f_{p,q}(k,\beta)$ is independent of $r$. In this case writing $x_r^{p,q} := (rp,rq)^{\Uk}$ to emphasise the dependence on $p$ and $q$, on the one hand the random variables $(W^{p,q}_r)_{r \geq 1}$ given by
\begin{align*}
W^{p,q}_r := \frac{1}{r} X_{0,r} = \frac{1}{ r p} \log Z_{x_0^{p,q} \to x_r^{p,q} }( \beta) 
\end{align*}
converge almost surely to $f_{p,q}(k,\beta)$ as $r \to \infty$. Now suppose $p'$ and $q'$ are another pair of integers such that $q'/p' = c$. Then it is immediate that $p'q = q'p$, and hence for every $r$
\begin{align*}
W_{p'r}^{p,q} = W_{pr}^{p',q'},
\end{align*}
and in particular, almost-surely we have
\begin{align} \label{aslim}
f_{p,q}(k,\beta) = \lim_{r \to \infty} W^{p,q}_{p'r} =  \lim_{s \to \infty}W^{p',q'}_{pr} = f_{p',q'}(k,\beta),
\end{align}
This ensures that $f_{p,q}(k,\beta)$ is independent of the choice of $p$ and $q$, and justifies us hereafter writing $f_c(k,\beta)$ for this quantity.
\ep

In Lemma \ref{special points}, for rational values of $c$ we defined $f_c(k,\beta)$ as a limit associated with the partition functions at a slope of $c$. The main remaining step in the proof of Theorem \ref{finite k thm} is the following lemma, which provides upper and lower bounds on the asymptotic growths of $k$-points at a $(1,c')$ asymptotics---where $c'$ is any positive number---in terms of those associated with diagonal points at rational tilts $c$.

\begin{lem} \label{key equations}

Let $\vx$ and $\vy$ be nice $k$-points, and let $c$ and $c'$ be positive reals such that $0 < c < c'$ and $c$ is rational. Then we have
\begin{align} \label{liminf}
\lim \inf_{ N \to \infty} \frac{1}{N} \log Z_{ \vx \to \vy + (N,c'N) }(\beta) \geq f_c(k,\beta) + k  (c' - c)  \beta \nu\qquad \text{almost-surely},
\end{align}
and 
\begin{align} \label{limsup}
\lim \sup_{ N \to \infty} \frac{1}{N} \log Z_{ \vx \to \vy + (N,c'N) }(\beta) \leq \frac{c'}{c} f_c(k,\beta) -  k  \left(\frac{c'}{c} - 1 \right)   \beta \nu \qquad \text{almost-surely}.
\end{align}
\end{lem}
\bp
First we prove equation \eqref{liminf}. Since the distribution of the random variable $\log Z_{ \vx \to \vy + (N,c'N) }(\beta) $ is invariant under translations $(\vx,\vy) \mapsto (\vx+  a, \vy + a)$, without loss of generality we may assume that $ 0^{ \Uk} \leq y$. Suppose $p$ and $q$ are positive integers such that $q/p = c$ and recall the notation
\begin{align*}
v_r := (rp,rq)^{\Uk}.
\end{align*} 
Since $\vx$ is nice, by Proposition \ref{diagonal bounding} there is a sufficiently large integer $r \in \mathbb{Z}$ such that $v_{r} \geq \vx$. Now for each $N > 0$, we may write
\begin{align*}
N = \alpha_N p + \beta_N,
\end{align*}
where $\alpha_N$ and $\beta_N$ are non-negative integers and $\beta_N \in \{0,1,\ldots,p-1\}$. Now consider the chain of nice $k$-points
\begin{align} \label{chain}
\vx \leq v_r \leq v_{\alpha_N} \leq (N, \alpha_N q )^{\Uk } \leq (N , c'N )^{\Uk} \leq \vy + (N, c'N),
\end{align}
where the fact that $(N , c'N )^{\Uk} \leq \vy + (N, c'N)$ is a consequence of translating the nequality $(0,0)^{\Uk} \leq \vy$. By applying the series bound to \eqref{chain} we have
\begin{align} \label{partition chain}
\frac{1}{N} \log Z_{ \vx \to \vy+  (N,c'N) }(\beta) \geq \frac{1}{N} \left( a + b_N + c_N + d_N + e_N\right)
\end{align}
where
\begin{align*}
a := \log Z_{ \vx \to \vy } (\beta), \qquad b_N := \log Z_{v_r \to v_{\alpha_N}} (\beta),  \qquad c_N := \log Z_{ v_{\alpha_N} \to (N,\alpha_N q)^{\Uk} }(\beta) \\
d_N := \log Z_{ (N,\alpha_Nq)^{\Uk} \to (N,c'N)^{\Uk}}(\beta),  \qquad e_N :=  \log Z_{ (N, c'N)^{\Uk} \to \vy + (N, c'N)}  (\beta).
\end{align*}
We now consider the asymptotic contribution coming from each of these five terms. First, we note that $a$ is independent of $N$, and hence 
\begin{align} \label{aeq}
\lim_{N \to \infty} \frac{1}{N} a_N = 0 \qquad \text{almost-surely.}
\end{align}
Second, we note that we may write $b_N := - p X_{r,\alpha_N}$, where $X_{r,s}$ is defined as Lemma \ref{special points}. In particular, by Lemma \ref{special points}, as $N \to \infty$, 
\begin{align} \label{beq}
\lim_{N \to \infty} \frac{1}{N} b_N = - \lim_{N \to \infty} \frac{1}{ \alpha_N p + \beta_N } p X_{r,\alpha_N } = - \lim_{N \to \infty} \frac{1}{\alpha_N} X_{r,\alpha_N } = f_c(k,\beta) \qquad \text{almost-surely.}
\end{align}
Now we consider the term $c_N$. Recalling that $v_{\alpha_N} := (\alpha_N p, \alpha_N q)^{\Uk}$, we note that $v_{\alpha_N}$ and $(N,\alpha_N q)^{\Uk}$ have the same vertical coordinates. It follows that there is only one $k$-path in $\Gamma_{ v_{\alpha_N} \to (N,\alpha_N q)^{\Uk}  }$ --- the $k$-path travelling horizontally. Moreover, each constituent $1$-path in the $k$-path has length $N - \alpha_N p \in \{0,1,\ldots,p-1\}$. In particular, for each $N$, $c_N$ is simply a sum of at most $k (p-1)$ independent and identically distributed random variables. It follows that 
\begin{align} \label{ceq}
\lim_{N \to \infty} \frac{1}{N} c_N = 0 \qquad \text{almost-surely.}
\end{align}
Turning our attention to $d_N$, we note that $ (N,\alpha_Nq)^{\Uk} $ and $(N,c'N)^{\Uk}$ have the same horizontal coordinates. In particular, there is only one $k$-path travelling between them, and hence $d_N$ is identical in distribution to the sum of $k \left( \lfloor c' N \rfloor -  \alpha_N q\right)$ independent random variables identically distributed like $\beta \omega(0,0)$. In particular, since
$
\lim_{N \to \infty} \frac{1}{N} k \left(  \lfloor c' N \rfloor -  \alpha_N q \right) = k (c' - c),
$
by applying the law of large numbers we have
\begin{align} \label{deq}
\lim_{N \to \infty} \frac{1}{N} d_N = k (c' - c) \beta \nu \qquad \text{almost-surely.}
\end{align}
Finally, by translation, the law of $e_N$ is independent of $N$. In particular,
\begin{align} \label{eeq}
\lim_{N \to \infty} \frac{1}{N} e_N = 0 \qquad \text{almost-surely.}
\end{align}
By plugging \eqref{aeq}, \eqref{beq}, \eqref{ceq}, \eqref{deq}, and \eqref{eeq},  into \eqref{partition chain}, we obtain the first equation, \eqref{liminf}.

\vspace{8mm}
Now we turn to proving the second equation, \eqref{limsup}. Since the proof of \eqref{limsup} is similar to the proof of \eqref{liminf}, we furnish fewer details in this case. Letting $\vx$ and $\vy$ be nice $k$-points, we assume without loss of generality that $\vy \leq 0^{\Uk}$. Moreover, there is a sufficiently small integer $r \in \mathbb{Z}$ such that $v_r \leq \vx$. 

For real numbers $u$, write $\lceil u \rceil$ for the least integer larger than $u$. Now for each integer $N$, define $\gamma_N := \lceil c' N/q \rceil$. Then
\begin{align} \label{floordef}
\lfloor c' N \rfloor = \gamma_N q - \delta_N,
\end{align}
where $\delta_N \in \{ 0,1,\ldots,q-1\}$. In particular $\lfloor c' N \rfloor \leq \gamma_N q$, and moreover $p \gamma_N =  p \lceil c' N/q \rceil \geq p \lceil c N / q \rceil = p \lceil N/ p \rceil \geq N$. Using the definition of $v_r$ to obtain the first inequality below, the fact that $\vy \leq 0^{\Uk}$ to obtain the second, $\lfloor c' N \rfloor \leq \gamma_N q$ to  obtain the third, and $p \gamma_N  \geq N$ in conjunction with Proposition \eqref{diagonal ordering} to obtain the fourth, it follows that
\begin{align} \label{chain2}
v_r \leq \vx \leq \vy + (N, c'N) \leq (N, c'N)^{\Uk} \leq ( \gamma_N p, c'N )^{\Uk} \leq v_{ \gamma_N}.
\end{align}
In particular, by applying the series bound to the chain of $k$-points \eqref{chain2} and rearranging, for every $N$ we have the inequality
\begin{align} \label{halloween gambit}
\frac{1}{N} \log Z_{ \vx \to \vy + (N, c'N)}(\beta) \leq \frac{1}{N} a_N - \frac{1}{N} b - \frac{1}{N} c_N - \frac{1}{N}d_N - \frac{1}{N} e_N,
\end{align} 
where 
\begin{align*}
a_N := \log Z_{ v_r  \to v_{\gamma_N}} (\beta), \qquad 
b := \log Z_{ v_r \to \vx}(\beta), \qquad
c_N := \log Z_{ y+ (N, c'N) \to (N, cN)^{\Uk} }(\beta),\\
d_N := \log Z_{ (N, c'N)^{\Uk} \to  ( \gamma_N p, c'N )^{\Uk} } (\beta), \qquad
e_N := \log Z_{  ( \gamma_N p, c'N )^{\Uk} \to v_{ \gamma_N} } (\beta).
\end{align*}
Now, in analogy with the proof of equation \eqref{liminf}, it is straightforward to show that
\begin{align} \label{tom}
\lim_{N \to \infty} \frac{1}{N} b = \lim_{N \to \infty} \frac{1}{N} c_N = \lim_{N \to \infty} \frac{1}{N} e_N = 0 \qquad \text{almost-surely.}
\end{align}
Now we consider the term $a_N$. In the context of Lemma \ref{special points}, we may write $a_N :=  - p X_{r , \gamma_N}$. Using \eqref{floordef} to obtain the second equality below, and Lemma \ref{special points} to obtain the third, we have
\begin{align} \label{dick}
\lim_{N \to \infty} \frac{1}{N} a_N := - \lim_{N \to \infty} \frac{p }{ N } X_{ r, \gamma_N } = - \frac{ c'}{c} \lim_{N \to \infty} \frac{1}{ \gamma_N} X_{r , \gamma_N } = \frac{c'}{c} f_c(k,\beta) \qquad \text{almost-surely.}
\end{align}
Finally, we consider the $d_N$ term. Since the vertical coordinates of $(N, c'N)^{\Uk}$ and $ ( \gamma_N p, c'N )^{\Uk}$ are identical, there is only a single $k$-path from $(N, c'N)^{\Uk}$ to $( \gamma_N p, c'N )^{\Uk}$, and in particular, $d_N$ is simply a sum of $k (p \gamma_N - N)$ independent and identically distributed random variables with the same law as $\beta \omega(0,0)$. Using the fact that $\lim_{N \to \infty}\frac{1}{N} k (p \gamma_N - N) = k \left( \frac{c'}{c} - 1 \right)$ and the law of large numbers, we have
\begin{align} \label{harry}
\lim_{N \to \infty} \frac{1}{N} d_N = k \left( \frac{c'}{ c} - 1 \right) \beta \nu \qquad \text{almost-surely.}
\end{align}
Combining \eqref{tom}, \eqref{dick} and \eqref{harry} in \eqref{halloween gambit}, we obtain \eqref{limsup}, completing the proof of Lemma \ref{key equations}.
\ep

We are now ready to complete the proof of Theorem \ref{finite k thm}.

\bp[Proof of Theorem \ref{finite k thm}]
First we recall that in Lemma \ref{special points}, for rational values of $c$ we defined a function $f_c(k,\beta)$ as the limit of a rescaled partition function involving diagonal $k$-points at the corners of the $(p,q)$ lattice, where $q/p = c$. We now show that this limit function is continuous in rational values of $c$, and deduce that $f_c(k,\beta)$ must also act as a limit function in a broader context.

First of all, if we assume $c'$ is also rational, and set $\vx = \vy = (0,0)^{\Uk}$, it is an immediate consequence of Lemma \ref{special points} and Lemma \ref{key equations} that for any fixed $k$ and $\beta$, $F(c) := f_c(k,\beta)$ satisfies
\begin{align} \label{ineq1}
F(c) + k(c'-c) \nu \beta \leq F(c') \leq \frac{c'}{c} F(c) - k \left( \frac{c'}{ c} - 1 \right) \beta \nu, \qquad \text{for $c < c'$ both rational}.
\end{align}
By rearranging \eqref{ineq1} we also have 
\begin{align} \label{ineq2}
\frac{c}{c'} F(c') + \lambda \left( 1 - \frac{c}{ c'} \right) \leq F(c) \leq F(c') - \lambda( c' - c), \qquad \text{for $c < c'$ both rational}.
\end{align}
In particular, \eqref{ineq1} and \eqref{ineq2} together imply that the function $f_c(k,\beta)$, so far only defined for rational values of $c$, is continuous in $c$. In particular, for \emph{any} $c' > 0$, we may define
\begin{align} \label{clim}
f_{c'}(k,\beta) := \lim_{ c \uparrow c' } f_c(k,\beta),
\end{align}
where the limit $c \uparrow c'$ in \eqref{clim} is taken through rational values $c$ tending up to $c'$. 

 We are now ready to show that the conclusions of Theorem \ref{finite k thm} hold for all nice $k$-points $\vx$ and $\vy$, and for all positive directions $c'$, with $f_{c'}(k,\beta)$ defined in \eqref{clim} acting as the limit. Indeed, by Lemma \ref{key equations}, on the one hand we have
\begin{align}  \label{aslim}
\lim \inf_{ N \to \infty} \frac{1}{N} \log Z_{ \vx \to \vy + (N,c'N) }(\beta) \geq \lim_{ c \uparrow c'} \left \{  f_c(k,\beta) + k  (c' - c)  \beta  \nu \right\} = f_{c'}(k ,\beta)\qquad \text{almost-surely},
\end{align}
and on the other hand
\begin{align}  \label{aslim2}
\lim \sup_{ N \to \infty} \frac{1}{N} \log Z_{ \vx \to \vy + (N,c'N) }(\beta) \leq \lim_{ c \uparrow c'} \left\{ \frac{c'}{c} f_c(k,\beta) -  k  \left(\frac{c'}{c} - 1 \right)   \beta \nu  \right\}  = f _{c'}(k, \beta) \qquad \text{almost-surely}.
\end{align}
(Both of the limits taken in \eqref{aslim} and \eqref{aslim2} are taken in rational values $c$ tending up to $c'$.) By combining \eqref{aslim} and \eqref{aslim2}, we have completed the proof of Theorem \ref{finite k thm}, which states that 
\begin{align*}
\lim_{ N \to \infty}\frac{1}{N} \log Z_{ \vx \to \vy + (N,c'N) }(\beta) = f_{c'}(k,\beta) \qquad \text{almost-surely.}
\end{align*}
\ep

In fact, from our proof of Theorem \ref{finite k thm}, we may immediately deduce Theorem \ref{ineq thm}:
\bp[Proof of Theorem \ref{ineq thm}]
Since the limit function $f_c(k,\beta)$ is continuous in $c$, the inequalities \eqref{ineq1} and \eqref{ineq2} hold for all pairs of positive numbers $c < c'$. In particular, by setting $k = 1$ we prove Theorem \ref{ineq thm}. 
\ep

\subsection{Proof of Theorem \ref{finite k thm 2}} \label{5.2}
We now prove Theorem \ref{finite k thm 2}, which states that $f_c(k,\beta) = k f_c(\beta),$ where $f_c(\beta) := f_c(1,\beta)$. In this section we will use the notation 
\begin{align*}
Z \left[ \vx \to \vy \right] := Z_{\vx \to \vy}( \beta)
\end{align*}
 since the $k$-points $\vx$ and $\vy$ appearing in here can be notationally heavy, and we will not be considering different values of $\beta$. The proof of Theorem \ref{finite k thm 2} is split into two inequalities, the first of which is significantly easier than the other. 

\bp[Proof of the inequality $f_c(k,\beta) \leq k f_c(\beta)$]
Consider setting $\vx = (0,0)^{\Uk}=  \vy$ in Theorem \ref{finite k thm}. Then by definition 
\begin{align} \label{apple1}
\lim_{N \to \infty} \frac{1}{N}  \E \log Z \left[ (0,0)^{\Uk} \to  (N,cN)^{\Uk} \right] = f_c(k,\beta)
\end{align}
On the other hand, the partition function for the analogous $1$-path satisfies
\begin{align} \label{apple2}
\lim_{N \to \infty}\frac{1}{N}   \E  \log Z \left[ (0,0) \to  (N,cN) \right] = f_c(1,\beta) =: f_c(\beta).
\end{align}
However, by the iterated parallel inequality \eqref{seperating}, for every pair of positive integers $m$ and $n$ we have the bound
\begin{align*}
 Z\left[  (0,0)^{\Uparrow k} \to (n,m)^{\Uparrow k} \right] \leq \prod_{i=1}^k Z \left[ (-i+1,i-1) \to (m-i+1,n+i-1) \right].
\end{align*}
Since each $Z \left[ (-i+1,i-1) \to (m-i+1,n+i-1) \right]$ has the same distribution as $Z \left[ (0,0) \to (n,m) \right]$, comparing \eqref{apple1} and \eqref{apple2} and  letting $N \to \infty$ with $m = cN$ yields the result.
\ep

\bp[Proof of the inequality $f_c(k,\beta) \geq k f_c( \beta)$]
Since both $f_c(k,\beta)$ and $f_c(\beta) := f_c(1,\beta)$ are continuous functions of $c$, it suffices to prove the result for rational $c$. To this end, let $p$ and $q$ be positive integers and let $c = q/p$. Then letting $\vx = \vy = (0,0)^{\Uk}$ in Theorem \ref{finite k thm} we have
\begin{align*}
\lim_{ M \to \infty} \frac{1}{M^2 p} \E[ \log Z \left[ (0,0)^{\Uk} \to (M^2 p, M^2 q )^{\Uk} \right] = f_c(k,\beta).
\end{align*}
The main idea of this proof is to show that the set of $k$-paths from $(0,0)^{\Uk}$ to $(M^2p, M^2 q )^{\Uk}$ contains a subset of $k$-paths travelling an distance order $M$ apart whose contribution is asymptotically close to $k f_c(\beta)$. More specifically, consider the $1$-points 
\begin{align*}
z_i^j(M) := (i-1) (0, Mq + 1) + j (Mp , Mq) \quad 1 \leq i \leq k, 1 \leq j \leq M - k ,
\end{align*}
and for $1 \leq j \leq M-k$ define associated $k$ points $\vz^j(M) := (z^j_1(M),\ldots,z^j_k(M))$. It is straightforward to check that 
\begin{align} \label{chainz}
(0,0)^{ \Uk } \leq \vz^{0}(M)\leq \ldots \leq \vz^{j-1}(M) \leq \vz^j(M) \leq \ldots \leq \vz^{ M - k }(M) \leq (M^2 p , M^2 q)^{\Uk}. 
\end{align}
We remark that since the points $z_i^j(M)$ form a lattice, the partition functions $Z \left[ \vz^{j-1}(M) \to \vz^j(M) \right]$ are identically distributed for different values of $j$. In particular, by applying the series bound \eqref{series bound} to partition functions associated with the chain of $k$-points in \eqref{chainz}, taking logarithms and then expectations, and then using the fact that $Z \left[ \vz^{j-1}(M) \to \vz^j(M) \right]$ are identically distributed, we have
\begin{align} 
&\frac{1}{M^2 p} \E[ \log Z \left[ (0,0)^{\Uk} \to (M^2 p, M^2 q )^{\Uk} \right]  \nonumber \\
 &\leq  \frac{1}{M^2p} \E \left[\log  Z \left[ (0,0)^{\Uk} \to \vz^{0}(M)\right] \right]  + \frac{1}{M^2p }  \E \left[ \log Z \left[ \vz^{M-k}(M)\to  (M^2p,M^2q)^{\Uk} \right] \right] \nonumber \\
 &+ \frac{M-k}{M^2p } \E \left[ \log Z \left[ \vz^{0}(M) \to \vz^{1}(M) \right] \right] \label{gluck},
\end{align}
We now consider the asymptotics of each of the three terms on the right hand side of \eqref{gluck}. First we look to use Corollary \ref{bluntbound} to obtain a lower bound for the quantity 
\[\lim \inf_{M \to \infty}   \frac{1}{M^2 p} \E \left[\log  Z \left[ (0,0)^{\Uk} \to \vz^{0}(M)\right] \right].\]
It is easy to see there is a sufficiently large constant $C_{p,q}$ not depending on $M$ such that $\varpi \left((0,0)^{\Uk}, \vz^{0}(M)  \right)  $, the total number of points contained in any $k$-path from $ (0,0)^{\Uk}$ to $ \vz^{0}(M)$, satisfies
\begin{align*}
\varpi \left((0,0)^{\Uk}, \vz^{0}(M)  \right)  \leq C_{p,q} M .
\end{align*} 
Now using Corollary \ref{bluntbound} to obtain the first inequality below, we have 
\begin{align*}
 \E \left[\log  Z \left[ (0,0)^{\Uk} \to \vz^{0}(M)\right] \right]   \geq \varpi  \left((0,0)^{\Uk}, \vz^{0}(M)  \right)  \beta \nu \geq - | \beta \nu | C_{p,q} M.
\end{align*} 
In particular, 
\begin{align} \label{glucka}
\lim \inf_{M \to \infty}   \frac{1}{M^2 p} \E \left[\log  Z \left[ (0,0)^{\Uk} \to \vz^{0}(M)\right] \right] \geq  \lim \inf_{M \to \infty}   \frac{1}{M^2} \left( - | \beta \nu | C_{p,q} M \right) = 0.
\end{align}
A similar argument allows us to obtain the bound for the term involving the partition functions of $k$-paths $\vz^{M-k}(M)$ to $(M^2p,M^2q)^{\Uk} $, giving us 
\begin{align} \label{gluckb}
\lim \inf_{M \to \infty}   \frac{1}{M^2 p} \E \left[\log  Z \left[ \vz^{M-k}(M)\to  (M^2p,M^2q)^{\Uk}  \right] \right] \geq 0.
\end{align}
Finally, we now consider the term involving $\E \left[ \log Z \left[ \vz^{0}(M)\to \vz^{1}(M) \right] \right]$. Suppose for each $1 \leq i  \leq k$, we have a collection $\pi_1,\ldots,\pi_k$ of $1$-paths such that each $\pi_i$ travels from $z^0_i(M) \to z^1_i(M)$. Note that for any $i$, the vertical coordinate of $z_{i+1}^0$ is given by $ i Mq + i$, which is strictly greater that the vertical coordinate of $z_i^1$, which is given by $i M q + i - 1$. This implies that every $k$-tuple of $k$-paths from $\vz^0(M)  \to \vz^1 (M)$ is guaranteed to be non-intersecting, and hence 
\begin{align*}
Z  \left[ \vz^{0}(M)\to \vz^{1}(M) \right]  = \prod_{ i =1}^k Z \left[ z^0_i (M) \to z^1_i (M) \right],
\end{align*}
and in particular,
\begin{align} \label{gluckc}
\E \left[ \log Z \left[ \vz^{0}(M)\to \vz^{1}(M) \right] \right] = k \E \left[ \log Z \left[ (0,0) \to (Mp , Mq) \right] \right].
\end{align}
Now using \eqref{glucka}, \eqref{gluckb} and \eqref{gluckc} in \eqref{gluck}, we have 
\begin{align} \label{sunny}
\lim \inf_{ M \to \infty} \frac{1}{M^2 p} \E[ \log Z \left[ (0,0)^{\Uk} \to (M^2 p, M^2 q )^{\Uk} \right] \geq  \lim \inf_{ M \to \infty} \frac{ k (M-k)}{ M^2 p }  \E \left[ \log Z \left[ (0,0) \to (Mp , Mq) \right] \right].
\end{align}
Finally, by using the case $k=1$ of Theorem \ref{finite k thm}, we have
\begin{align} \label{winter}
\lim_{ M \to \infty} \frac{1}{M p}  \E \left[ \log Z \left[ (0,0) \to (Mp , Mq) \right] \right] = f_c(\beta).
\end{align}
Plugging \eqref{winter} into \eqref{sunny}, we yield
\begin{align*}
\lim \inf_{ M \to \infty}  \frac{1}{M^2 p} \E[ \log Z \left[ (0,0)^{\Uk} \to (M^2 p, M^2 q )^{\Uk} \right] \geq k f_c(\beta),
\end{align*} 
which completes the proof of the inequality $f_c(k,\beta) \geq k f_c(\beta)$.

\ep

\subsection{Proof of Theorem \ref{many k} } \label{5.3}
We now prove Theorem \ref{many k} by a direct application of Kingman's subadditive ergodic theorem, Theorem \ref{kin thm}.
\bp[Proof of Theorem \ref{many k}]
Let $h = (h^1,h^2)$ such that $h^1 \leq 0 < h^2$. For $i \geq 0$, define
\begin{align*}
x_i := (i-1)h, ~~ y_i = y + (i-1)h.
\end{align*}
For positive integers $r < s$, define the $(s-r)$-points $\vx^{r,s} $ and $\vy^{r,s}$ by 
\begin{align*}
\vx^{r,s} := (x_r,x_{r+1},\ldots,x_{s-1}), ~~ \vy^{r,s} = (y_r,x_{r+1},\ldots,y_{s-1} ),
\end{align*}
and define the random variable
\begin{align*}
X_{r,s} := \log Z \left[ \vx^{r,s} \to \vy^{r,s} \right].
\end{align*}
We now verify that with the collection of random variables $(X_{r,s})_{0 \leq r < s }$ are in the set up of Theorem \ref{kin thm}, Kingman's subadditive ergodic theorem. Note first that by the parallel bound \eqref{parallel bound}, for any $ r < s < t$ we have the subadditivity
\begin{align*}
X_{r,s} + X_{s,t } \geq X_{r,t},
\end{align*}
which proves the first condition. The second condition is immediate since the weights are independent and identically distributed. Finally, the third point followings from Corollary \ref{bluntbound}.  To see this, it is simple to verify that $\vx^{0,r} \leq \vy^{0,r}$ and $\varpi(\vx^{0,r}, \vy^{0,r}) = r || y- x||_1$, and hence by Corollary \ref{bluntbound} 
\begin{align*}
\E [ X_{0,r} ] =  \E \left[  \log Z \left[ \vx^{0,r} \to \vy^{0,r} \right] \right] \geq \beta \nu || y - x||_1 r.
\end{align*}
It follows that we are in the setting of Theorem \ref{kin thm}, and hence
\begin{align*}
\frac{1}{k} X_{0,k} = \frac{1}{k} \log Z \left[ \vec{x}^{h \uparrow k} \to \vy^{h \uparrow k}  \right]
\end{align*}
converges to a deterministic limit $J_{x,y,h}(\beta)$ depending on $x$, $y$, $h$, and $\beta$. Since the weights of the polymer are independent and identically distributed, it is plain that $J$ only depends on $x$ and $y$ through the difference  $y-x$. This completes the proof of Theorem \ref{many k}.
\ep

In the sequel we return to writing $Z_{\vx \to \vy}(\beta)$ (as opposed to $Z[\vx \to \vy]$) for the $\vx$ to $\vy$ partition function.

\subsection{Proof of Theorem \ref{toeplitz thm}} \label{toe proof}

To prove Theorem \ref{toeplitz thm}, we recall Szeg{\"o}'s limit theorem, following B{\"o}ttcher and Silbermann \cite[Chapter 5]{BS} and Bump \cite[Chapter 42]{bum}. Szeg{\"o}'s limit theorem is concerned with the asymptotics of the determinants of Toeplitz matrices --- those matrices of the form of the form $(d_{j-i})_{1 \leq i,j \leq k}$. 

To set up the limit theorem, for every collection of coefficients $(d_m)_{m \in \mathbb{Z}}$ we associate a Laurent series $f(s) = \sum_{m \in \mathbb{Z}} d_m s^m$ representing a function $f: \mathbb{T} \to \mathbb{C}$ on the unit circle. This Laurent series is known as the symbol of the Toeplitz matrix.

We define the Wiener and Besov norms of the symbol $f$ by 
\begin{align}
|| f ||_{W} := \sum_{ m \in \mathbb{Z}} |d_m| < \infty ~~~\text{and}~~~ || f ||_{B_2^{1/2}} := \sum_{m \in \mathbb{Z}} (|m|+1)|d_m|^2 < \infty
\end{align}
respectively.
We say $f \in W \cap B_2^{1/2}$ if both $|| f ||_{W} $ and $|| f ||_{B_2^{1/2}} $ are finite.  Recalling the definition of winding number stated in \eqref{wind def}, we now state Szeg{\"o}'s limit theorem.

\begin{thm}[Szeg{\"o}'s limit theorem] \label{szego thm}
Suppose $f(s) := \sum_{m \in \mathbb{Z}} d_m s^m$ is a Laurent series in $W \cap B_2^{1/2}$ with no zeroes on $\mathbb{T}$ and winding number zero. Moreover, suppose there is a Laurent series $\sum_{m \in \mathbb{Z}} c_m s^m$ satisfying
\begin{align*}
f(s) = \exp \left(  \sum_{m \in \mathbb{Z}} c_m s^m \right).
\end{align*}
Then 
\begin{align*}
\lim_{ k \to \infty} \left(  e^{ -k  c_0} \det_{i,j=1}^k \left( d_{j-i} \right)  \right) = \exp \left( \sum_{m=1}^\infty m c_m c_{-m} \right). 
\end{align*}
\end{thm}

We now use the Lindstr{\"o}m-Gessel-Viennot formula to relate the asymptotics of many-path infinite temperature polymer to the asymptotics of Toeplitz determinants. 
Now, for $k$-points $\vx$ and $\vy$, by the \eqref{gv} 
\begin{align} \label{toe?}
Z_{ \vx \to \vy}(0 ) = \det_{i,j=1}^k (Z_{x_i \to y_j}(0)).
\end{align}
We are interested in circumstances where the matrix  $\left( Z_{x_i \to y_j}(0)\right)_{1 \leq i,j \leq k}$ is a Toeplitz matrix for every $k$, which occurs whenever each $Z_{x_i \to y_j}(0)$ is a function of $j-i$ only, and the $k$-points $\vx$ and $\vy$ have the form
\begin{align} \label{shifted}
\vx_i  = x + (i-1)h, ~~ \vy_i = y + (i-1), ~~ i = 1,\ldots,k.
\end{align}
for some base points $x$ and $y$ in $\mathbb{Z}^2$, and some direction $h \in \mathbb{Z}^2$. That is, $\vx$ and $\vy$ are of the form $\vx^{h\uk}$ and $\vy^{h\uk}$.

Indeed, when $\vx = \vx^{h\uk}$ and $\vy = \vy^{h\uk}$, each $Z_{x_i \to y_j}(0) =  d_{j-i}$, where 
\begin{align*}
d_m := \binom{ z^1 + z^2+ m(h^1 + h^2 ) }{ z^1 + m h^1 },
\end{align*}
where $z := y-x$, and it follows that the partition function may be written as a Toeplitz determinant: $Z_{ \vx \to \vy} ( 0 ) = \det_{i,j=1}^k ( d_{j-i})$.  Now define the symbol
\begin{align*}
a_{z,h}(s) := \sum_{ m \in \mathbb{Z}} d_m s^m. 
\end{align*}
Whenever $h_1 < 0 < h_2$, it is clear that only finitely many of the $d_m$ are non-zero, and hence $a_{z,h}(s)  := \sum_{ m \in \mathbb{Z}} d_m s^m$ is contained in $W \cap B_S^{1/2}$. By the assumption in Theorem \ref{toeplitz thm}, $a_{z,h}$ has winding number $0$, and it folows 
that we are in the set up of the Szeg{\"o}'s limit theorem, Theorem \ref{szego thm}. In particular, if $(c_m)_{m \in \mathbb{Z}}$ are such that
\begin{align*}
a_{z,h}(s) = \exp \left( \sum_{ m \in \mathbb{Z}} c_m s^m \right),
\end{align*}
then
\begin{align*}
\lim_{ k \to \infty} \left(  e^{ -k  c_0} Z_{\vx^{h \uparrow k} \to \vy^{h \uparrow k} }(0) \right) = \lim_{ k \to \infty} \left(  e^{ -k  c_0} \det_{i,j=1}^k ( d_{j-i})\right) = \exp \left( \sum_{m=1}^\infty m c_m c_{-m} \right),
\end{align*}
completing the proof of Theorem \ref{toeplitz thm}.

\subsection{The asymptotics of the high-temperature scaled-$k$ limit} \label{infinite temperature}
This section is devoted to proving Lemma \ref{infinite growth}. In this direction, first we state --- and, for completeness, include a proof of --- the following result regarding the infinite temperature partition function on a finite rectangle. We note that, via the well-known bijection between non-intersecting paths and plane partitions, this is equivalent to Macmahon's formula. We refer the reader to Stanley \cite[Section 7.20]{sta} for an approach using the RSK correspondence.

\begin{prop} \label{comb lem}
Let $m,n,k$ be positive integers satisfying $k \leq m \wedge n$. Then
\begin{align} \label{final}
Z_{ (1,1)^{\uk} \to (n,m)^{\downarrow k} }(0)  =\frac{ H(m + n - k) H(k) H(m-k) H(n-k) }{ H(m) H(n) H(m+n - 2k) },
\end{align}
where $H(N):= \prod_{j=0}^{N-1} j!$ is the superfactorial.
\end{prop}

\bp
Let $m,n,k$ be positive intgers satisfying $k \leq m \wedge n$, and suppose $\pi = (\pi_1,\ldots,\pi_k)$ is a $k$-path going from $(1,1)^{\uk} \to (n,m)^{\downarrow k}$. It is easily seen that for each $i$, the $i^{\text{th}}$ path $\pi_i$ must go through the points $x_i = (k+1-i,i)$ and $y_i := (m-i+1,n-k+i)$. Let $\vx = (x_1,\ldots,x_k)$ and $y = (y_1,\ldots,y_k)$ be the associated $k$-points. Then we have the following identity for the rectangular zero-temperature partition functon
\begin{align} \label{path q}
Z_{ (1,1)^{\uk} \to (n,m)^{\downarrow k} }(0) = Z_{ \vx \to \vy }(0). 
\end{align}
By the Lindstr{\"o}m-Gessel-Viennot formula \eqref{gv}, 
\begin{align*}
Z_{ \vx \to \vy }(0)  = \det_{i,j = 1}^k \left( Z_{x_i \to y_j} (0) \right).
\end{align*}
Now note that the partition function $Z_{x_i \to y_j}(0)$ simply counts the number of paths from $x_i$ to $y_j$, that is
\begin{align} \label{path q 2}
Z_{x_i \to y_j} (0) = \binom{ m + n - 2k  }{ n - k + j - i }.
\end{align}
It bares remarking at this stage that the right-hand-side of \eqref{path q} can be expanded as a more tractable determinant than the left-hand-side, since the path length from each $x_i$ to each $y_i$ is $m + n - 2k$, which does not depend on $i$ or $j$.\\

According to a well known determinant identity (see Krattenthaler \cite[Equation (2.17)]{kra}), 
\begin{align} \label{kra eq}
\det_{i,j = 1}^k \left( \binom{a+b}{a + j- i} \right) = \prod_{r = 1}^k \prod_{s=1}^a \prod_{t=1}^b \frac{r+s+t - 1}{ r + s + t - 2}.
\end{align}
Moreover, it easily verified that
\begin{align} \label{barnes id}
\prod_{r = 1}^k \prod_{s=1}^a \prod_{t=1}^b \frac{r+s+t - 1}{ r + s + t - 2} = \frac{ H(k + a + b) H(k) H(a) H(b) }{ H(k+a) H(k+b) H(a+ b) },
\end{align} 
where $H(N) := \prod_{j=1}^{N-1}  j!$ is the superfactorial. Setting $a = m -k$ and $b = n-k$, and combining \eqref{path q}, \eqref{path q 2}, \eqref{kra eq} and \eqref{barnes id}, we obtain \eqref{final}.
\ep

\bp[Proof of Lemma \ref{infinite growth}]
To prove Lemma \ref{infinite growth}, it remains to study the large-$N$ asymptotics of \eqref{final} when $m = cN$, $n=N$ and $k=\alpha N$, for any $c > 0$ and $\alpha \leq c \wedge 1$. By adapting \eqref{bar asy}, we see that for $p > 0$,
\begin{align} \label{HpN}
\log H( p N) = N^2 \left( \frac{p^2}{2} \log p + \frac{p^2}{2} \log N - \frac{3 p^2 }{4} + o(1)\right). 
\end{align}
Now the result follows from \eqref{final} and \eqref{HpN}, noting the identity
\begin{align*}
(c+1 - \alpha)^2 + \alpha^2 + (c - \alpha)^2 + (1- \alpha)^2 - c^2 - 1^2 - (c+1-2\alpha)^2 = 0.
\end{align*}
\ep

\section{Proofs of results in Section \ref{SIM}} \label{SIM proofs}
\subsection{Proofs of Proposition \ref{interface 2} and Proposition \ref{eigen interface}}  \label{interface proofs}

In this section we prove Propositions \ref{interface 2} and \ref{eigen interface}, which state that certain random functions related to random polymers and random matrices may be expressed as stochastic interfaces.

\bp[Proof of Proposition \ref{interface 2}]
Using \cite[Equations (2.8) and (2.9)]{COSZ}, taking the change of variable $y_i = \log u_i$ in \cite[Equation (1.3)]{COSZ}, and using the definition of Whittaker functions in \cite[Section 3.2]{COSZ} we find that the diagonal $\varphi_\mu^N$ has law 
\begin{align*}
 \frac{1}{\Gamma(\mu)^{N^2}} \exp \left(  - e^{ - \lambda _N } - \mu \sum_{i =1}^N \lambda_i \right) g^{\mathsf{exp}}(\lambda)^2 \prod_{ i =1}^N d \lambda_i,
\end{align*}
which by Lemma \ref{integrate out}, agrees with the law of the corresponding stochastic interface model. It remains to show that the processes have the same law not just on the diagonal but everywhere on $S_N$. In this direction, by \cite[Theorem 3.7 (ii)]{COSZ}, the conditional law of $\{ \varphi_\mu^N(i.j) : 1 \leq i \leq j \leq N \}$ given the values on the diagonal $\left( \varphi_\mu^N(1,1),\ldots,\varphi_\mu^N(N,N) \right)= (\lambda_1,\ldots,\lambda_N)$ is given by 
\begin{align*}
\frac{1}{ g^{\mathsf{exp}}(\lambda) } \exp \left( - \sum_{ \langle x,y \rangle } e^{\phi(y) - \phi(x)} \right)   \prod_{ 1 \leq i < j \leq N } d \phi(i,j) \prod_{i=1}^N \delta_{\lambda_i}(d \phi(i,i)).
\end{align*}
Finally, it follows from the construction in \cite{COSZ} that the random variables $\{ \tau^N(m,k) \}_{ 1 \leq k \leq m \leq N } $ are conditionally independent of the random variables $\{ \tilde{\tau}^N(m,k) \}_{ 1 \leq k \leq m \leq N }$ given the diagonal. By the distributional symmetry of $\varphi_\mu^N$ on either side of the diagonal, this completes the proof. 
\ep

\bp[Proof of Proposition \ref{eigen interface}]
We prove the result for the eigenvalue process associated with the Gaussian unitary ensemble, omitting a proof for the Laguerre unitary ensemble, which is almost identical. 

Let $\varphi^N:S_N \to \mathbb{R}$ be a stochastic interface with interaction potential $\mathsf{bead}$ and weight functions $W^{\mathsf{GUE}}_{i,j}(u) = \ind_{i=j} \frac{1}{2} u^2$. By Lemma \ref{integrate out} and the continuus Gelfand-Tsetlin volume formula \eqref{wdf}, it is plain that the diagonal $(\lambda_1,\ldots,\lambda_N) := (\varphi^N(1,1),\ldots,\varphi^N(N,N))$ is distributed according to the probability measure
\begin{align} \label{gin}
\frac{1}{Z H(N)^2 } \Delta_N(\lambda)^2 \exp \left( - \frac{1}{2} \sum_{ i = 1}^N \lambda_i^2 \right) \prod_{ i = 1}^N d \lambda_i.
\end{align}
Equation \eqref{gin} is known as the Ginibre formula, and it is well known that the $N$ eigenvalues $\lambda_1 > \ldots > \lambda_N$ of $H$ of a matrix from the Gaussian Unitary Ensemble also have this law (see e.g. Mehta \cite[Theorem 3.3.1]{meh}). In other words, the stochastic interfaces $\varphi^N$ and $\varphi^N_{H,U}$ have the same distribution on the diagonal $i=j$.

It is now a consequence of Baryshnikov \cite[Theorem 0.7]{bar} that the equality in distribution for $\varphi^N_{H,U}$ and $\varphi^N$ holds not just on the diagonal but everywhere on the square $S_N$. 
\ep

\subsection{Proof of Theorem \ref{large mu}} \label{large mu proof}
In this section we consider the large-$\mu$ asymptotics of the interface $\varphi_\mu^N$, providing a proof of Theorem \ref{large mu}.

\bp[Proof of Theorem \ref{large mu}]

By definition, the law of the interface $\varphi^N_\mu:S_N \to \mathbb{R}$ is proportional to
\begin{align*}
\exp \left( - e^{ - \phi(N,N) } - \mu \sum_{i=1}^N \phi(i,i) - \sum_{ \langle x, y \rangle \in S_N^* } e^{ \phi(y) - \phi(x) } \right).
\end{align*}
We want to take a change of variables so that the interaction term competes with the weight term for large $\mu$. To this end, consider the change of variables
\begin{align} \label{psi def 2}
\theta^N_\mu(i,j) := \varphi^N_\mu(i,j)  +  ( 2N + 1 - j - i) \log \mu.
\end{align}
It is straightforward to show that the random function $\theta^N_\mu:S_N \to \mathbb{R}$ is itself a stochastic interface whose law is proportional to $\exp( - \mu \mathcal{F}^N[\theta])$, where 
\begin{align} \label{cal f 2}
\mathcal{F}^N[ \phi ] :=  e^{ - \phi(N,N)} + \sum_{i=1}^N \phi(i,i) + \sum_{ \langle x, y \rangle \in S_N^* } e^{ \phi(y) - \phi(x) }.
\end{align}
Plainly, as $\mu \to \infty$ the interface $\theta^N_\mu$ converges in distribution to the deterministic function $\theta^N_{\mathsf{min}}:S_N \to \mathbb{R}$ minimising $F^N[\theta]$. 
That completes the proof that $\theta_\mu^N$ converges in distribution to a deterministic function $\theta_{\mathsf{min}}^N$ minimising a variational problem. It remains to show that this minimiser has the form \eqref{psi form}.

Since the law of $\theta_\mu^N$ is symmetric about the diagonal $i=j$, the limit function $\theta_{\mathsf{min}}$ must also be symmetric, and without loss of generality we consider $\theta_{\mathsf{min}}(i,j)$ for $i \leq j$. 

To compute $\theta_{\mathsf{min}}(i,j)$, we consider the direct implications for the polymer partition functions of sending $\mu \to \infty$. For integer values of $\mu$, $\frac{1}{\zeta_\mu(z)}$ can be written as a sum of $\mu$ independent and identically distributed exponential variables with mean $1$. It follows that as $\mu \to \infty$
\begin{align} \label{mu conv}
\frac{1}{ \mu \zeta_\mu(z)} \text{ converges in distribution to $1$}.
\end{align}
Since each product in the sum $\tau^N_\mu(m,k):= \sum_{ \pi \in \Gamma^N(m,k) } \prod_{ z \in \pi} \zeta_\mu(z)$ contains $k(N+m-k)$ weights, we have
\begin{align} \label{tau conv}
\mu^{k(N+m-k)} \tau^N_\mu(m,k) \text{ converges in distribution to $ \# \Gamma^N(m,k) $}.
\end{align}
Now using the definitions of $\varphi_\mu^N$ and $\theta_\mu^N$  (which appear in \eqref{phi def} and \eqref{psi def 2} respectively), for $i \leq j$ we have
\begin{align} \label{xi exp}
\theta^{N}_{ \mu}(i,j) &:= \varphi^N_\mu(i,j) + (  2N  + 1 - (j + i))\log \mu \nonumber \\ 
 &:= \log \left( \frac{\mu^{ i(N+(N-j+i)-i)}~ \tau^N_\mu(N-j+i,i)  }{ \mu^{(i-1)(N + (N-j+i) - (i-1) ) } ~ \tau^N_\mu(N-j+i,i-1)} \right).
\end{align}
Combining \eqref{xi exp} with \eqref{tau conv}, it follows that as $\mu \to \infty$
\begin{align} 
\theta^N_\mu(i,j) ~\text{converges in distribution to $\log \left( \frac{\# \Gamma^N(N-j+i,i)  }{\# \Gamma^N(N-j+i,i-1) } \right)$. }
\end{align}
It remains to compute
\begin{align} \label{comb def}
\theta^N_{\mathsf{min}}(i,j) = \log \left( \frac{\# \Gamma^N(N-j+i,i)  }{\# \Gamma^N(N-j+i,i-1) } \right).
\end{align}
In Lemma \ref{comb lem}, we gave the following explicit expression for the cardinality of each $ \Gamma^N(m,k)$
\begin{align} \label{final 2}
\# \Gamma^N(m,k)   =\frac{ H(N + m -k) H(k) H(N - k) H(m-k ) }{ H(N) H(m) H(N+m-2k ) }.
\end{align}
Combining \eqref{comb def} with \eqref{final 2}, and using the definition $H(n) := \prod_{i=1}^N (i-1)!$ of the superfactorial, after some calculation we obtain the formula
\begin{align} \label{psi form 2}
\theta^N_{\mathsf{min}}(i,j) := \log \left( \frac{ (i-1)! (2N - j - i + 1)! (2N - j - i)! }{(2N-j)! (N-j)! (N-i)! } \right), ~~~ i \leq j,
\end{align}
completing the proof.
\ep

\section{Proofs of results in Section \ref{thermodynamics} } \label{thermodynamics proofs}

\subsection{Proof of Theorem \ref{big N big mu thm}} \label{big N big mu proof}
In this section we prove Theorem \ref{big N big mu thm}.

\bp[Proof of Theorem \ref{big N big mu thm}]
Since $\theta^N_{\mathsf{min}}$ is symmetric, the limit function of the rescaled interface $\bar{\theta}^N_{\mathsf{min}}$ must also be symmetric, and hereafter without loss of generality we consider $\{ s \leq t \}$. In this case, by \eqref{psi form} and Definition \ref{rescaled}, the rescaled interface associated with $\theta_{\mathsf{min}}^N$ is given by 
\begin{align*}
\bar{\theta}^N_{\mathsf{min}}(s,t) = \frac{1}{N} \log \left( \frac{ (i-1)! (2N - j - i + 1)! (2N - j - i)! }{(2N-j)! (N-j)! (N-i)! } \right) ~~ \text{for $ i = \lf sN \rf$ and $j = \lf t N \rf$.}
\end{align*}
By Stirling's formula, $\frac{1}{N} \log ( \lf u N \rf ) =u (\log N - 1) + q(u) + o(1)$, where $q(u) := u \log u$. Using the fact that
\begin{align*}
s + 2(2-s-t) - (2-t) - (1-t) - (1-s) = 0,
\end{align*}
it follows that $\bar{\theta}^N_{\mathsf{min}}(s,t) \to \xi_{\mathsf{ht}}(s,t) $, given by
\begin{align*}
\xi_{\mathsf{ht}}(s,t) = q(s) + 2q(2-s-t) - q(2-t) - q(1-t) - q(1-s).
\end{align*}

\ep

\subsection{Calculations surrounding Conjecture \ref{conj gauss}} \label{gauss proof}

Section \ref{gauss proof} is dedicated to sketching calculations leading us to Conjecture \ref{conj gauss}, which anticipaties Gaussian fluctuations of the logarithmic partition functions of the log-gamma polymer at high temperature. 

First of all, consider that by the central limit theorem, for each $z$, as $\mu \to \infty$ the random variable
\begin{align*}
N_\mu(z) := \frac{ 1/\zeta_\mu(z) - \mu }{\sqrt{\mu}}
\end{align*}
converges in distribution to a standard normal random variable. Now define the variable $r_z(\mu)$ through the equation
\begin{align} \label{r def}
1 + \mu^{-1/2} r_z(\mu) = \mu \zeta_\mu(z).
\end{align}
It is straightforward to check that $r_z(\mu) = \left(  1 + \mu^{-1/2}N_\mu(z) \right)^{-1} N_z(\mu) $, and hence $r_\mu(z)$ also converges in distribution to a standard normal random variable as $\mu \to \infty$. 

Assuming \eqref{lim dens} and \eqref{dec}, we now consider the asymptotics of the partiton functions under a suitable scaling as $\mu$ and $N$ tend to infinity. Consider the random variable
\begin{align*}
G^N_\mu(m,k) &:= \frac{ \mu^{ k(N+m-k)} \tau^N_\mu (m,k)}{ \# \Gamma^N(m,k) } 
\end{align*}
Using \eqref{r def} and the fact that each path in $\Gamma^N(m,k)$ contains $k(N+m-k)$ weights, we have
\begin{align*}
G^N_\mu(m,k)  &= \frac{1}{ \# \Gamma^N(m,k) } \sum_{ \pi \in \Gamma^N(m,k) } \prod_{ z \in \pi} \left( 1 + \frac{1}{ \sqrt{\mu} }  r_z(\mu) \right)\\
&= \frac{1}{ \# \Gamma^N(m,k) } \sum_{ \pi \in \Gamma^N(m,k) } \sum_{F \subset \pi} \prod_{ z \in F}  \frac{1}{ \sqrt{\mu} }  r_z(\mu)\\
&= \sum_{F \subset S_N} \mu^{ - \# F/2} Q^N_{m,k}( F \in \pi)  \prod_{ z \in F}  r_z(\mu),
\end{align*}
where the sum in the final line above is taken over all subsets $F$ of the square $S_N := \{1,\ldots,N\}^2$. Expanding further, and considering the $n!$ ways of ordering the elements of a set $F$ of size $n$, we have
\begin{align} \label{G eq}
G_\mu^N(m,k) = \sum_{n = 0}^\infty \frac{ \mu^{-n/2} }{n!} \sum_{ (z_1,\ldots,z_n) \in S_N^{(n)} } Q_{m,k}^N( \{z_1,\ldots,z_n\} \in \pi  ) \prod_{i=1}^n r_\mu( z_i),
\end{align}
where the internal sum is taken over all distinct $n$-tuples of elements of $S_N$. When $N$ is large, for each $n$ the set of non-distinct $k$-tuples $S_N^n - S_N^{(n)}$ is small compared to the set of distinct $k$-tuples $S_N^{(n)}$. Using this fact as well as the asymptotic decoupling \eqref{dec}, we anticipate that when $N$ is large, and both $m$ and $k$ have order $N$, then
\begin{align*}
 \sum_{ (z_1,\ldots,z_n) \in S_N^{(n)} } Q_{m,k}^N( \{z_1,\ldots,z_n\} \in \pi  ) \prod_{i=1}^n r_\mu( z_i) = (1 + o(1)) \left( \sum_{ z \in S_N }Q_{m,k}^N( z \in \pi  ) r_\mu(z) \right)^n,
\end{align*}
in which case using \eqref{G eq} reduces to 
\begin{align} \label{lim def}
G^N_\mu(m,k)   = (1 + o(1)) \sum_{n=0}^\infty \frac{1}{n!}  \left(\mu^{-1/2} \sum_{ z \in S_N}Q_{m,k}^N( z \in \pi  ) r_\mu(z) \right)^n.
\end{align} 
Taking the scaling $\mu = \kappa N^2$, and using \eqref{lim def}, we anticipate that 
\begin{align} \label{whitenoise}
\log G^N_{\kappa N^2} ( cN, \alpha N) \to \frac{1}{ \kappa}  \int_S w_{c, \alpha} ( u,v) \dot{W}(u,v),
\end{align}
where $\dot{W}$ is a space-time white noise on $S$ independent of $c$ and $\alpha$. (We remark that by considering $\tilde{\tau}$ in place of $\tau$, we also expect 
\begin{align*}
\log \tilde{G}^N_{\kappa N^2} (cN, \alpha N ) := \log \left(\frac{ \mu^{ k(N+m-k)} \tilde{\tau}^N_\mu (m,k)}{ \# \tilde{\Gamma}^N(m,k) } \right)
\end{align*}
converges in distribution to $\frac{1}{ \kappa} \int_S \tilde{w}_{c , \alpha} (u,v) \dot{W}(u,v)$, where $\tilde{w}$ is defined in analogy to \eqref{lim dens}, instead using the uniform measure $\tilde{Q}^N_{m,k}$ on $\tilde{\Gamma}^N(m,k)$. )

Now define 
\begin{align} \label{H def}
H_{\kappa}^N(s,t) := 
\begin{cases}
\log G^N_{\kappa N^2} ( tN, s N) ~~ s \leq t.\\
\log \tilde{G}^N_{\kappa N^2} (sN, t N )  ~~ s \geq t. \\
\end{cases}
\end{align}
Then by \eqref{whitenoise},  we anticipate that $H_\kappa^N$ converges in distribution to a Gaussian process $H_\kappa(s,t)$ on $S$ with convariance
\begin{align*}
\E \left[  H_\kappa(s,t) H_\kappa(s',t') \right ] = \frac{1}{ \kappa^2} \big \langle h_{s,t} , h_{s',t'} \big \rangle_{ L^2(S)} 
\end{align*}
where for $s \leq t$ let $h_{s,t}(u,v):= w_{t,s}(u,v)$, and define $h_{t,s}(u,v) := h_{s,t}(v,u)$, and $\langle f ,g \rangle_{ L^2(S)} := \int_S f(z) g(z) dz$.

\section{Derivation of the surface tension of the bead model} 
\label{stproof}

\subsection{Coordinate changes and Wulff functionals}

It will be useful to define the change of coordinates mapping the set $T := \{ ( x, y) : 0 \leq x \leq y \leq 1 \}$ to $U := \{ (r,\tau) : 0 \leq \tau \leq 1, ~ \tau/2 \leq r \leq 1-\tau/2 \}$ defined by setting
\begin{align}\label{tilted coor}
\tau := y-x ~ \text{and} ~r := \frac{ y + x }{ 2}.
\end{align}
We say that $(r,\tau)$ are the tilted coordinates. Note that the transformation $(x,y) \mapsto (r,\tau)$ has unit Jacobian, and for suitable functions $f$ we have
\begin{align} \label{tilted diff}
\frac{ \partial f}{ \partial x } = \frac{1}{2}  \frac{ \partial f}{ \partial r } -  \frac{ \partial f}{ \partial \tau } \qquad \text{and} \qquad \frac{ \partial f}{ \partial y } = \frac{1}{2}  \frac{ \partial f}{ \partial r } +  \frac{ \partial f}{ \partial \tau }.
\end{align}
It turns out to be useful to study the surface tension in the tilted change of coordinates $(r, \tau)$. In light of \eqref{tilted diff}, we  define the surface tension associated with an interaction potential $V$ in the tilted coordinates by
\begin{align} \label{to tilted}
\sigma^V_{\mathsf{tilted}}( p,q) := \sigma^V  \left(  \frac{1}{2} p - q, \frac{1}{2}p + q \right).
\end{align}
Now with these definitions at hand, we define the triangular Wulff functional $\mathcal{W} : \mathcal{C}^1(U) \to \mathbb{R}$ associated with the interaction potential $V$ by 
\begin{align}  \label{wulff funct}
\mathcal{W}^V[ v ] := \int_{U}  \sigma^V_{ \mathsf{tilted}} \left( \frac{ \partial  v}{ \partial r}, \frac{ \partial v}{ \partial \tau} \right) d r d \tau.
\end{align}
Given a function defined on the bottom of the triangle $U$, we would like to study the minimum energy attained by a shape defined on $U$ agreeing with this function on the boundary. Namely, define the minimal Wulff functional $\mathcal{M}^V:\mathcal{C}^1([0,1]) \to \mathbb{R}$ by
\begin{align} \label{min wulff}
\mathcal{M}^V[\rho] := \min_{ v \in \mathcal{C}^1(U)_\rho  } \mathcal{W}^V[ v]  =  \min_{v \in \mathcal{C}^1(U)_{ \rho} } \int_{U} \sigma^V_{\mathsf{tilted}} \left( \frac{ \partial v}{ \partial r} , \frac{ \partial v}{ \partial \tau} \right) dr d \tau,
\end{align}
where the minimisation is taken over $\mathcal{C}^1(U)_\rho := \{ v \in \mathcal{C}^1(U) : v(0,r) = \rho(r) ~\forall~ r \in [0,1] \}$. 

We will be interested in functions $\rho: [0,1] \to \mathbb{R}$ minimising $\mathcal{M}^V[\rho]$ subject to certain costs. (Of course, here there are two different steps of minimisation taken place, first over internal values on the triangle $T$ or $U$, and then on the boundary shape in $[0,1]$.) In any case, in order to study these sort of problems, we need a functional calculus allowing us to differentiate functionals. To this end, let $A:\mathcal{C}^1([0,1]) \to \mathbb{R}$ be a functional, and define the functional derivative $D_A$ of $A$ to be the map $D_A: \mathcal{C}^1([0,1]) \times \mathcal{C}^1([0,1]) \to \mathbb{R}$ given by 
\begin{align*}
D_A( \rho, \eta) :=  \frac{d}{du} A[\rho + u \eta] \Big|_{u = 0}.
\end{align*}
It is easily verified that $D_A$ is linear in the second argument for every $\rho$, and hence by the Riesz representation theorem \cite{rie}, for every $\rho$ there exists a Radon measure $\Lambda_A(\rho, \cdot)$ on $[0,1]$ such that
\begin{align*}
D_A( \rho, \eta) = \int_0^1 \eta(r) \Lambda_A(\rho, dr).
\end{align*}
We call the measure $\Lambda_A(\rho,\cdot)$ the Riesz measure associated with $A$ and $\rho$. Note that whenever $\rho$ is a local minima or maxima of the functional $A$, the Riesz measure $\Lambda_A(\rho,\cdot)$ is zero. Let us also remark that whenever $A:\mathcal{C}^1([0,1]) \to \mathbb{R}$ is itself a linear map we have
\begin{align} \label{linear linear}
D_A(\rho,\eta) = A[\eta].
\end{align}
Now to see an example of the functional derivative, consider the thermodynamic Vandermonde determinant 
\begin{align*}
\bm{\Delta}[ \rho] := \int_{0 < s < t < 1 } \log( \rho(s) - \rho(t) ),
\end{align*}
which was introduced in Section \ref{stbead}. A simple calculation taking the derivative inside the integral verifies that the Riesz measure associated with $\bm{\Delta}$  has continuous density given by 
\begin{align} \label{riesz van}
\Lambda_{\bm{\Delta}}(\rho, r) = \int_0^1  \frac{ 1 }{ \rho(r) - \rho(s) } ds.
\end{align}
(That is, for measurable $\eta$ we have $\int \eta(r) \Lambda_{\bm{\Delta}}(\rho, r)  dr = \int \eta(r) \int_0^1  \frac{ 1 }{ \rho(r) - \rho(s) } ds dr$.) 
We will be interested in studying the functional derivatives of minimal Wulff functionals. Namely, given a surface tension $V$, and a minimal Wulff functional $\mathcal{M}^V$, suppose that for a given $\rho$ the Riesz measure $\Lambda_{\mathcal{M}^V} \left( \rho, \cdot \right)$ has a continuous density which we denote by $\Lambda(\rho,r)$ for $r \in [0,1]$. 

We now show that given $\rho$, the shape $v^*: U \to \mathbb{R}$ minimising the Wulff functional $\mathcal{W}^V$ subject to $v^*(r, 0) = \rho(r)$ must satisfying the differential equation
\begin{align} \label{wulff eq}
\frac{ \partial \sigma^V_{\mathsf{tilted}}}{ \partial q} \left( \rho'(r), \frac{ \partial v^*}{ \partial \tau} (0,r) \right) + \Lambda(\rho,r) = 0 \qquad \text{for $r \in [0,1]$.}
\end{align}

\bp[Derivation of \eqref{wulff eq}]
Given a function $\rho:[0,1] \to \mathbb{R}$, consider breaking the minimisation procedure into two steps, by writing
\begin{align} \label{zeroth min}
\min_{ v \in \mathcal{C}^1(U)_\rho } \mathcal{W}^V [ \rho] = \min_{ f  } \min_{ v \in \mathcal{C}^1(U)_{ \rho, f} } \mathcal{W}^V [ v ] ,
\end{align}
where the former minimisation on the right hand side is taken over continuously differentiable functions $f:[0,1] \to \mathbb{R}$, and the latter minimisation on the right hand side is taken over $  \min_{ v \in \mathcal{C}^1(U)_{ \rho, f} }$, the set of continuously differentiable $v:U \to \mathbb{R}$ satisfying 
\begin{align} \label{bcs}
v(r,0) = \rho(r) \qquad \text{and} \qquad \frac{ \partial v}{ \partial \tau} ( r, 0 ) = f(r). 
\end{align}
Now given functions $\rho, f$, consider finding a $v$ minimising $\mathcal{W}^V[ v]$ subject to the boundary conditions \eqref{bcs}. For subsets $I$ of $[0,1]$, let $U_I := U \cap \{ \tau \in I \}$ and consider the disjoint union $U = U_{[0,h) } \cup U_{[h,1]}$. Indeed, by seperately considering the integration over the regions in this disjoint union, for small $h$ we have 
\begin{align} \label{first min}
&  \min_{ v \in \mathcal{C}^1(U)_{ \rho, f} } \mathcal{W}^V [ v ] \nonumber \\
&=  o(h) + h \int_0^1 \sigma_{\mathsf{tilted}} \left( \rho'(r), f(r) \right)  dr + \min_{ v \in \mathcal{C}^1(U_{[h,1]})_{ \rho + h f}} \int_{ U_{[h,1]} }\sigma^V_{\mathsf{tilted}} \left( \frac{ \partial v}{ \partial r} , \frac{ \partial v}{ \partial \tau} \right) dr d \tau,
\end{align}
where $\mathcal{C}^1(U_{[h,1]})_{ g }$ denotes the set of continuously differentiable functions defined on $U_{[h,1]}$ and satisfying $v(h,r) = g(r)$ for all $- h/2 < r < 1 - h/2$. For small $h$, it is easy to check from geometric considerations that 
\begin{align} \label{geoh}
\min_{ v \in \mathcal{C}^1(U_{[h,1]})_{ g }} \int_{ U_{[h,1]} }\sigma^V_{\mathsf{tilted}} \left( \frac{ \partial v}{ \partial r} , \frac{ \partial v}{ \partial \tau} \right) dr d \tau = ( 1 - 2h ) \mathcal{M}^V[ g ] + o(h).
\end{align}
(This follows from the fact that $\Leb(U_{[h,1]}) = (1 - 2h + o(h)) \Leb( U)$, where $\Leb$ denotes the Lebesgue measure on $\mathbb{R}^2$.) 
In particular,
\begin{align}  \label{small h}
\min_{ v \in \mathcal{C}^1(U_{[h,1]})_{ \rho + h f}} \int_{ U_{[h,1]} } \sigma^V_{\mathsf{tilted}} \left( \frac{ \partial v}{ \partial r} , \frac{ \partial v}{ \partial \tau} \right) dr d \tau   = \mathcal{M}^V[\rho] + h \left( - 2 \mathcal{M}^V[\rho ] + D_{\mathcal{M}^V}(\rho,f) \right) + o(h).
\end{align}
Plugging \eqref{small h} into \eqref{first min}, we obtain
\begin{align} \label{second min}
 \min_{ v \in \mathcal{C}^1(U)_{ \rho, f} } \mathcal{W}^V [ v ] := \mathcal{M}^V[\rho] + h \left\{ \int_0^1 \sigma_{\mathsf{tilted}} \left( \rho'(r), f(r) \right)  dr   - 2 \mathcal{M}^V[\rho ] + D_{\mathcal{M}^V} (\rho,f) \right\} + o(h)
\end{align} 
Comparing \eqref{second min} and \eqref{zeroth min}, we see that the minimising function $f$ in 
\eqref{zeroth min} must minimise  $\mathcal{G}_\rho[f ]$, where
\begin{align*}
\mathcal{G}_\rho[f ] := \int_0^1 \sigma_{\mathsf{tilted}} \left( \rho'(r), f(r) \right)  dr  + D_{\mathcal{M}^V} (\rho,f).
\end{align*} 
Using our assumption that the Riesz derivative has a continuous density $\Lambda_{\mathcal{M}^V}(\rho, r) $, we may write
\begin{align*}
\mathcal{G}_\rho[f] := \int_0^1 \left\{    \sigma_{\mathsf{tilted}}\left( \rho'(r), f(r) \right)  + \Lambda_{\mathcal{M}^V}(\rho, r) f(r)  \right\} dr .
\end{align*}
The result now follows from a standard variational argument. Indeed, suppose $f^*$ is a minimiser of $\mathcal{G}_\rho[f]$. Then for every test function $\eta$ we have $\frac{ d}{ds} \mathcal{G}_\rho[ f^* + s \eta] = 0 $, which amounts to 
\begin{align} \label{var rep}
\int_0^1 \left\{ \frac{ \partial \sigma_{\mathsf{tilted}}}{ \partial q} \left( \rho'(r) , f^*(r) \right)  +  \Lambda_{\mathcal{M}^V}(\rho, r)  \right\} \eta(r) dr.
\end{align}
In order for \eqref{var rep} to hold for every $\eta$, we must have $\frac{ \partial \sigma_{\mathsf{tilted}}}{ \partial q} \left( \rho'(r) , f^*(r) \right)  +  \Lambda_{\mathcal{M}^V}(\rho, r) = 0$. This equation must now be satisfied by $f^*(r) := \frac{ \partial v}{ \partial \tau} (r,0)$, where $v^*$ is the minimiser of $\mathcal{W}^V[v]$ over $\mathcal{C}^1(U)_\rho$, completing the derivation of \eqref{wulff eq}.
\ep

\subsection{A scaling property of the surface tension}
The surface tension of the bead model satisfies a certain scaling property. Namely it is easily seen by setting $V(u) = \mathsf{bead}(u)$ in  \eqref{finite surface} that for all $s,t < 0$, and $\lambda > 0$, we have 
\begin{align} \label{finite bead scaling}
\sigma^{\mathsf{bead}}_N( \lambda s, \lambda t) = - \log \lambda + \sigma^{\mathsf{bead}}_N( s,t),
\end{align}
for each of the finite surface tensions  $\sigma^{\mathsf{bead}}_N$. The limit $\sigma^{\mathsf{bead}}(s,t) := \lim_{N \to \infty} \sigma^{\mathsf{bead}}_N(s,t)$ inherits the same property, and using \eqref{to tilted}, so does the tilted surface tension:
\begin{align} \label{bead scaling}
\sigma^{\mathsf{bead}}_{\mathsf{tilted}}( \lambda p ,\lambda q) = - \log \lambda + \sigma_{\mathsf{tilted}}^{\mathsf{bead}}( p,q).
\end{align}
In particular, for $p < 0$ by setting $\lambda = |p|$ in \eqref{bead scaling}, we see that there is a function $\Omega$ such that 
\begin{align} \label{omega rep}
\sigma^{\mathsf{bead}}_{\mathsf{tilted}} (p, q) = - \log (|p| )  +  \Omega(q/|p| ).
\end{align}
We remark that in particular,
\begin{align} \label{omega rep 2}
\frac{ \partial \sigma^V_{\mathsf{tilted}}}{ \partial q} (p, q)  = \frac{1}{|p|} \Omega'(q/|p|).
\end{align}

\subsection{Derivation of equation \ref{st}}
Recall now that the goal of this section is to derive the formula
\begin{align} \label{st2}
\sigma_{\mathsf{tilted}}^{\mathsf{bead}}( p, q) := 
\begin{cases}
 - \log \left( |p| \cos \left( \pi \frac{ q}{ |p| } \right) \right) & \text{if} ~ p < 0, |q| < |p| \\
\infty &\text{otherwise},\\
\end{cases}
\end{align}
for the surface tension of the bead model. First, we remark that whenever we do not have $p < 0$ and $q < |p|$, we have $\sigma_{\mathsf{tilted}}^{\mathsf{bead}}( p, q) := \infty$. Indeed, it is immediate from the definition that for the finite (untilted) surface tension we have $\sigma^{\mathsf{bead}}_N( s,t ) < \infty$ if and only if both $s$ and $t$ are strictly negative. In particular, $\sigma^{\mathsf{bead}}(s,t) := \lim_{N \to \infty} \sigma^{\mathsf{bead}}_N( s,t ) $ also inherits this property. Now  setting $s = \frac{1}{2} p - q$ and $t \frac{1}{2} p  + q$, we see that $(s < 0, t < 0)$ is equivalent to $(p < 0, |q| < |p|)$. In particular,
\begin{align} \label{finiteness}
\sigma_{\mathsf{tilted}}^{\mathsf{bead}}( p, q) < \infty \qquad \text{only if } p < 0, |q| < |p|.
\end{align}
It remains to actually compute $\sigma_{\mathsf{tilted}}^{\mathsf{bead}}( p, q)$ in the case that $ p < 0$ and $|q| < |p|$. 

Clearly, by \eqref{omega rep}, \eqref{omega rep 2} and \eqref{finiteness}, in order to establish \eqref{st2} it is remains to show that
\begin{align} \label{tanrep}
\Omega'(s) = \pi \tan \pi s \qquad \text{for $s \in (-1,1)$}.
\end{align}

To this end, recall that in Section \ref{stbead} we introduced the semicircle law as the minimiser of a variational problem. We now restate this problem here in terms of the tilted coordinates system. Namely given $x \in [0,1]$, we recall that $\rho_{\mathsf{sc}}(x)$ is the solution in $[-2,2]$ to the equation
\begin{align} \label{rhodef2}
\int_{\rho_{\mathsf{sc}}(x)}^2  \frac{1}{2\pi} \sqrt{ 4 - u^2 } du = x.
\end{align}
Now consider the function $\xi_{\mathsf{sc}}$ in \eqref{psi sc} defined in the tilted change of coordinates \eqref{tilted coor}, namely
\begin{align} \label{xidef2}
\xi_{\mathsf{sc}} (  r , \tau )  := \sqrt{1 - \tau} \rho_{\mathsf{sc}} \left( \frac{r - \tau/2}{ 1- \tau}  \right) .
\end{align}
Let $K := \{ (r,\tau) : \tau \in [-1,1] , |\tau|/2 < r < 1 - |\tau|/2 \}$ be the diamond in $\mathbb{R}^2$. Restating the variational characterisation of $\xi_{\mathsf{sc}}$ in Section \ref{stbead} here in tilted coordinates,  the function $\xi_{\mathsf{sc}}$ is the minimiser over $v:K \to \mathbb{R}$ of the functional
\begin{align*}
E[v] := \int_K \sigma^V_{\mathsf{tilted}} \left( \frac{ \partial v}{ \partial r} , \frac{ \partial v}{ \partial \tau} \right) dr d \tau + \frac{1}{2}  \int_0^1 v(r,0)^2 dr. 
\end{align*}
Equivalently, the restriction of $\xi_{\mathsf{sc}}$ to $U = \{ (r, \tau) : \tau \in [0,1], \tau/2 < r < 1 - \tau/2 \}$ is the minimiser of 
\begin{align*}
E'[v] := 2 \int_U \sigma^V_{\mathsf{tilted}} \left( \frac{ \partial v}{ \partial r} , \frac{ \partial v}{ \partial \tau} \right) dr d \tau + \frac{1}{2} \int_0^1 v(r,0)^2 dr. 
\end{align*}

Note that $\xi(r,0 ) = \rho_{\mathsf{sc}}(r)$. In particular, breaking the minimisation of $E'[v]$ into two steps, we see that the function $\xi(r,\tau)$ is the minimiser of $\int_U \sigma^V_{\mathsf{tilted}} \left( \frac{ \partial v}{ \partial r} , \frac{ \partial v}{ \partial \tau} \right) dr d \tau$ over $\mathcal{C}^1(U)_{\rho_{\mathsf{sc}}}$. It follows that by \eqref{wulff eq}, we must have 
\begin{align} \label{wulff eq 2}
\frac{ \partial \sigma^V_{\mathsf{tilted}}}{ \partial q} \left( \rho_{\mathsf{sc}}'(r), \frac{ \partial \xi_{\mathsf{sc}} }{ \partial \tau} (0,r) \right) + \Lambda^{\mathsf{bead}}(\rho_{\mathsf{sc}},r) = 0 \qquad \text{for every $r \in [0,1]$.}
\end{align}

 Using the fact that the equation \eqref{wulff eq 2} involves each of the known quantities $\rho_{\mathsf{sc}}$, given in \eqref{rhodef2}, $\xi_{\mathsf{sc}}$ given in \eqref{xidef2} and $\Lambda^{\mathsf{bead}}$ given in \eqref{riesz van}, along with the representation \eqref{omega rep}, this is sufficient to establish \eqref{tanrep}. 

Indeed, as a first step, using \eqref{omega rep 2} and \eqref{riesz van} we may rewrite \eqref{wulff eq 2} as
\begin{align} \label{big final}
\frac{1}{ \rho_{\mathsf{sc}}'(r) } \Omega' \left( \frac{  \frac{ \partial \xi_{\mathsf{sc}}}{ \partial \tau} (r,0) }{ \rho_{\mathsf{sc}}'(r)} \right) + \int_0^1 \frac{1}{ \rho_{\mathsf{sc}}(r) - \rho_{\mathsf{sc}}(s) } ds = 0  \qquad \text{for every $r \in [0,1]$.}
\end{align}

It turns out to be easiest to consider equation \eqref{big final} through the change of variable $\rho_{\mathsf{sc}}(r) = 2 \sin \phi$. Now some computations are required. Let $f_{\mathsf{sc}}(x) = \frac{1}{2\pi} \sqrt{ 4 - x^2 } \ind_{x \in [-2,2]}$, and let $F_{\mathsf{sc}}(x) := \int_{-\infty}^x f(u) du$. Then $\rho_{\mathsf{sc}}$ is the inverse function of $1 - F$. Straightforward geometric considerations tell us that for $\phi \in [-\pi/2,\pi/2]$, 
\begin{align*}
F_{\mathsf{sc}}( 2 \sin \phi) = \frac{1}{2} + \frac{ \phi}{ \pi } + \frac{ \sin \phi \cos \phi}{ \pi}.
\end{align*}
Using the inverse function theorem to obtain the first inequality below, it follows that with $\rho_{\mathsf{sc}}(r) = 2 \sin \phi$ we have
\begin{align} \label{bobo}
\rho_{\mathsf{sc}}'(r) = - \frac{1}{F_{\mathsf{sc}}'(\rho_{\mathsf{sc}}(r))} = - \frac{1}{ \frac{1}{2\pi} \sqrt{ 4 - (2 \sin \phi)^2 } } =-  \frac{ \pi }{ \cos \phi} .
\end{align}
A computation shows that
\begin{align*}
\frac{ \partial \xi_{\mathsf{sc}} }{ \partial \tau} (\tau,r) = \frac{ - 1/2}{ (1-\tau)^{1/2} } \rho_{\mathsf{sc}} \left( \frac{ r - \tau/2}{ 1- \tau} \right) + \left(  \frac{ r- \tau/2}{ (1-\tau)^{3/2}} -  \frac{ 1/2}{ (1-\tau)^{1/2} }\right)  \rho_{\mathsf{sc}}' \left( \frac{ r - \tau/2}{ 1- \tau} \right),
\end{align*}
and in particular, on the line $\tau = 0$ we have
\begin{align} \label{giraffe}
\frac{ \partial \xi_{\mathsf{sc}} }{ \partial \tau} (r, \tau) = - \frac{1}{2} \rho_{\mathsf{sc}}(r) +  \left(r - \frac{1}{2}\right) \rho'_{\mathsf{sc}}(r).
\end{align}
Now note that if $\rho_{\mathsf{sc}}(r) = 2 \sin \phi$, then $r = 1 -  F_{\mathsf{sc}}(\rho_{\mathsf{sc}}(r)) =\frac{1}{2} - \frac{ \phi}{ \pi } - \frac{ \sin \phi \cos \phi}{ \pi}$. Thus using \eqref{bobo}, \eqref{giraffe} reduces to
\begin{align} \label{bobo 2}
\frac{ \partial  \xi_{\mathsf{sc}}}{ \partial \tau}(r,\tau) = \frac{-1}{2}2 \sin \phi + \left( \frac{ \phi}{ \pi } + \frac{ \sin \phi \cos \phi}{ \pi} \right) \frac{ \pi }{ \cos \phi } = \frac{ \phi}{ \cos \phi}.
\end{align}
We are now ready to show that equation \eqref{big final} implies \eqref{tanrep}. On the one hand, by \eqref{bobo} and \eqref{bobo 2} with $\rho_{\mathsf{sc}}(r) = 2 \sin \phi$ we  have
\begin{align} \label{eX}
\frac{1}{ \rho_{\mathsf{sc}}'(r) } \Omega' \left( \frac{  \frac{ \partial \xi_{\mathsf{sc}}}{ \partial \tau} (r,0) }{ \rho_{\mathsf{sc}}'(r)} \right)  = -\frac{ \cos \phi }{ \pi} \Omega' \left( - \phi/ \pi\right).
\end{align}
On the other hand, we now claim that thanks a to special property of the semicircle law, we have 
\begin{align} \label{eY}
\int_0^1 \frac{1}{ \rho_{\mathsf{sc}}(r) - \rho_{\mathsf{sc}}(s) } ds = \frac{1}{2} \rho_{\mathsf{sc}}( r ) = \sin \phi.    
\end{align}
Indeed, by changing variable at the point $\rho_{\mathsf{sc}}(r) = 2 \sin \phi$, we have
\begin{align} \label{eY}
\int_0^1 \frac{1}{ \rho_{\mathsf{sc}}(r) - \rho_{\mathsf{sc}}(s) } ds = - \int_{-2}^2 \frac{1}{ 2 \sin \phi - u } \frac{1}{2\pi } \sqrt{ 4 - u^2 } du. 
\end{align}
It is now a property of the semircle law well known in random matrix theory (see e.g. \cite[Equation 4.2.7]{meh}) that for $\alpha \in [-2, 2]$ we have$\int_{-2}^2 \frac{1}{ \alpha - u } \frac{1}{2\pi } \sqrt{ 4 - u^2 } du = \frac{1}{2} \alpha$, and hence by \eqref{eY} we have 
\begin{align} \label{eZ}
\int_0^1 \frac{1}{ \rho_{\mathsf{sc}}(r) - \rho_{\mathsf{sc}}(s) } ds = - \sin \phi
\end{align}
(where $\phi$ is implicit in the equation $\rho_{\mathsf{sc}}(r) = 2 \sin \phi$).

Plugging \eqref{eX} and \eqref{eZ} into \eqref{big final}, we obtain
\begin{align*}
 \Omega'( - \phi/\pi) = - \pi \tan \pi \phi,
\end{align*} 
which implies that $\Omega'(s) = \pi \tan \pi s$, proving \eqref{tanrep}, and therefore the formula \eqref{st}.

\end{document}